\documentclass[11pt]{article}

\usepackage[latin1]{inputenc} 
\usepackage{amsthm,amsmath,amssymb,bm} 
\usepackage{graphicx}
\usepackage{color}
\usepackage{verbatim}

\usepackage{enumitem} 
\setlist[enumerate]{topsep=0pt,itemsep=-1ex,partopsep=1ex,parsep=1ex}

\theoremstyle{plain}
\newtheorem{theo}{Theorem}[section]

\newtheorem{lemma}[theo]{Lemma}

\theoremstyle{definition}
\newtheorem{defn}[theo]{Definition}

\newcommand{\mc}[1]{\mathcal{#1}}
\newcommand{\mb}[1]{\mathbb{#1}}
\newcommand{\nib}[1]{\noindent {\bf #1}}

\newcommand{\bsize}[1]{\left| #1 \right|}

\newcommand{\bgen}[1]{\left\langle #1 \right\rangle}
\newcommand{\sgen}[1]{\langle #1 \rangle}

\newcommand{\sub}{\subseteq}

\newcommand{\Ra}{\Rightarrow}
\newcommand{\sm}{\setminus}

\newcommand{\es}{\emptyset}
\newcommand{\pl}{\partial}

\newcommand{\car}{\circlearrowright}

\newcommand{\aA}{\alpha}
\newcommand{\bB}{\beta}
\newcommand{\gG}{\gamma}
\newcommand{\dD}{\delta}

\newcommand{\lL}{\lambda}
\newcommand{\tT}{\theta}
\newcommand{\sS}{\sigma}
\newcommand{\oO}{\omega}
\newcommand{\ups}{\upsilon}

\newcommand{\aB}{{\bm{\alpha}}}
\newcommand{\iB}{{\bm{i}}}

\newcommand{\GG}{\Gamma}

\newcommand{\OO}{\Omega}

\newcommand{\Ss}{\Sigma}

\newcommand{\LL}{\Lambda}

\def\qed{\hfill $\Box$}

\title{Coloured and directed designs}

\author{Peter Keevash\thanks{Mathematical Institute,
University of Oxford, Oxford, UK. Email: keevash@maths.ox.ac.uk.
\newline \hspace*{1.8em}Research supported
in part by ERC Consolidator Grant 647678.}}

\begin{document}

\maketitle

\begin{abstract}
We give some illustrative applications of our recent 
result on decompositions of labelled complexes,
including some new results on decompositions of
hypergraphs with coloured or directed edges.
For example, we give fairly general conditions
for decomposing an edge-coloured graph into rainbow triangles,
and for decomposing an $r$-digraph into tight $q$-cycles.
\end{abstract}

\begin{center}
{\em To L\'aszl\'o Lov\'asz on his seventieth birthday}
\end{center}

\section{Introduction}

When can we decompose an object into copies of some other object?
This vague question suggests a number of mathematical problems.
Within graph theory, a fundamental instance of this question
asks for a decomposition (i.e.\ partition of the edge set)
of the complete graph $K_n$ into copies of $K_q$.
We require $n \ge q^2-q+1$ by Fisher's inequality
(see e.g.\ \cite[Theorem 19.6]{vLW}).
If $q$ is one more than a prime power then the
lines of a projective plane give a construction
with $n=q^2-q+1$, but we do not know any construction
with $n=q^2-q+1$ when $q$ is not of this form;
the Prime Power Conjecture suggests that there are none.
On the other hand, we may fix $q$ and ask for conditions
on $n$ that guarantee a decomposition (perhaps only for 
large $n>n_0(q)$ so as to exclude the difficulties
associated with the Prime Power Conjecture).
The first such result, obtained by Kirkman in 1846 (see \cite{RobinW}),
shows that $K_n$ has a triangle decomposition iff
$n$ is $1$ or $3$ modulo $6$.

\begin{figure}
\includegraphics{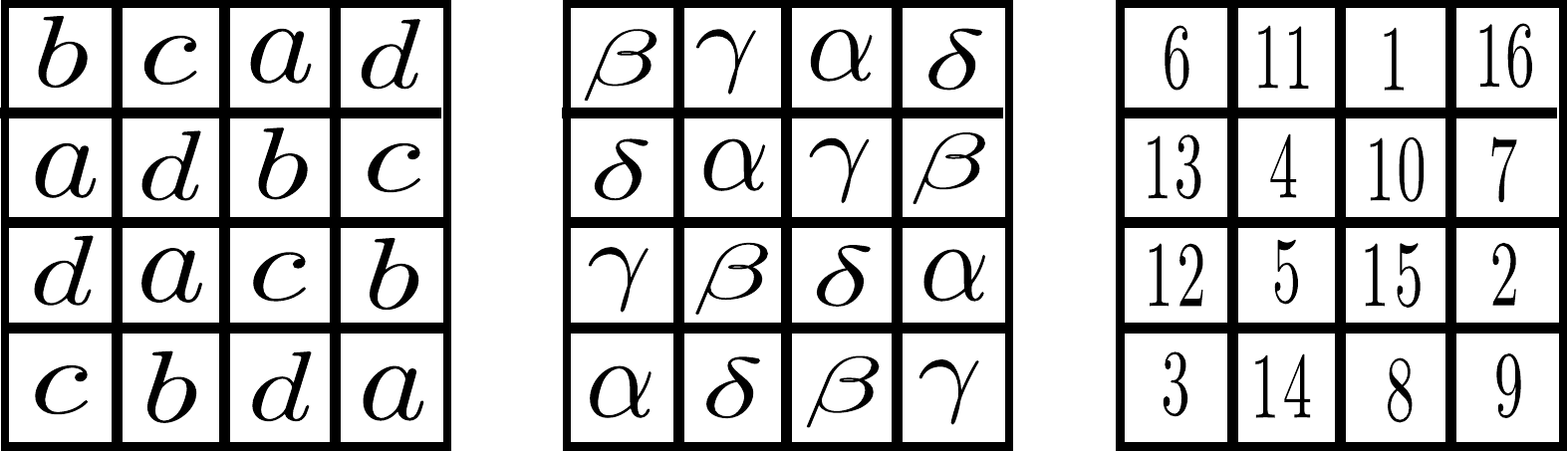}
\caption{Orthogonal and magic squares}
\label{fig:magic}
\end{figure}

These beginnings suggest several possible directions
for further generalisation. From the combinatorial perspective
(taken in this paper), one may ask for a decomposition of $G$
by copies of $H$ where $G$ and $H$ are any given graphs,
or hypergraphs, or indeed other related structures
(we will consider coloured and directed hypergraphs).
On the other hand, the above questions also have natural
interpretations in Design Theory, which suggests many further questions 
(some of which also have natural combinatorial interpretations).
Perhaps the oldest topic in this area is that of Latin and Magic squares,
which have their roots in antiquity (see \cite[Chapter 2]{CD});
they were given prominence in the Western mathematical tradition 
by Euler in 1776, who posed the {\em 36 officer's puzzle}, 
which was open until its solution by Tarry in 1900.
In modern terminology, the result is that there is no
pair of orthogonal Latin squares of order 6.
A pair of orthogonal Latin squares of order 4 is
illustrated in Figure \ref{fig:magic}, together 
with an associated magic square (obtained by assigning
values $1,2,3,4$ to $a,b,c,d$ and $0,4,8,12$ to $\aA,\bB,\gG,\dD$).

\begin{figure}
\begin{center}
\includegraphics{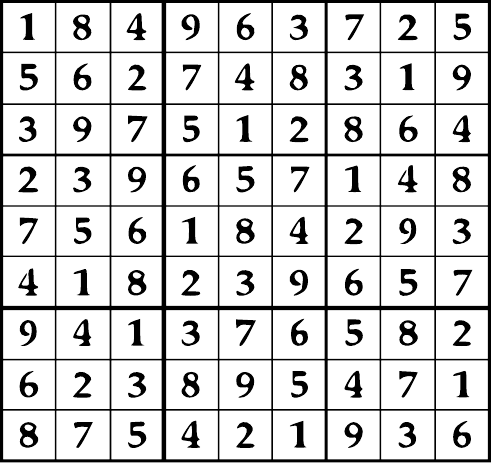}
\end{center}
\vspace{-0.3cm}
\caption{A completed Sudoku puzzle}
\label{fig:sudoku}
\end{figure}

In general, a Latin square of order $n$ is a labelling of the cells
of an $n$ by $n$ square with $n$ symbols so that every symbol
appears once in each row and once in each column.
An equivalent combinatorial description
is a triangle decomposition of $K_3(n)$, the complete tripartite graph
with parts of size $n$. Indeed, we identify the three parts with the sets
of rows, columns and symbols of the square, and then each cell corresponds
to a triangle in the obvious way. For a pair of orthogonal Latin squares 
of order $n$ we require two such squares with the extra condition that
every pair of symbols appears together once; 
this is analogously equivalent to a $K_4$-decomposition of $K_4(n)$
(and similarly for larger numbers of mutually orthogonal Latin squares).
We have chosen the pair in Figure \ref{fig:magic} with the extra property
that both diagonals use all symbols in both squares, so as to obtain
a magic square (all rows, columns and diagonals have the same sum). 
In Figure \ref{fig:sudoku} we illustrate the popular puzzle
of completing a partially filled Sudoku square,
which is a Latin square of order $9$ partitioned into $3$ by $3$
subsquares each of which uses every symbol once.

We now consider the generalisations of the above problems from graphs 
to $r$-graphs (hypergraphs in which every edge has size $r$).
When does an $r$-multigraph $G$ have a decomposition into copies
of some fixed $r$-graph $H$? The case that $H=K^r_q$ is the 
complete $r$-graph on $q$ vertices is of particular interest, as
a $K^r_q$-decomposition of $K^r_n$ is equivalent to a Steiner $(n,q,r)$ system,
i.e.\ a collection of blocks of size $q$ in a set of size $n$
covering every set of size $r$ exactly once. For example,
if $(q,r)=(3,2)$ a triangle decomposition of $K_n$ 
is equivalent to a Steiner Triple System. More generally,
giving each edge of $K^r_n$ some fixed multiplicity $\lL$,
a $K^r_q$-decomposition of $\lL K^r_n$ is equivalent to 
a $(n,q,r,\lL)$ design. Some necessary conditions for the existence
of a $K^r_q$-decomposition of an $r$-multigraph $G$ may be observed 
by considering the degrees. The degree of $e \sub V(G)$ is the number of
edges of $G$ containing $e$, i.e.\ the size of the neighbourhood
$G(e) = \{ f \sub V(G) \sm e: e \cup f \in G \}$.
We say $G$ is {\em $K^r_q$-divisible} if $|G(e)|$ is divisible
by $\tbinom{q-|e|}{r-|e|}$ for all $e \sub V(G)$;
this is a necessary condition for a $K^r_q$-decomposition,
as every copy of $K^r_q$ containing $e$ contains 
$\tbinom{q-|e|}{r-|e|}$ edges that contain $e$.
For example, a necessary condition for the existence of 
a $(n,q,r,\lL)$ design is $\tbinom{q-i}{r-i} \mid \lL \tbinom{n-i}{r-i}$
for all $0 \le i \le r-1$. The Existence Conjecture, proved in \cite{Kexist},
is that if $n>n_0(q,r,\lL)$ is large and this divisibility condition holds
then there is a $(n,q,r,\lL)$ design. More generally, we can find a
$K^r_q$-decomposition in any $K^r_q$-divisible $r$-multigraph $G$ that
is sufficiently dense and quasirandom.

The Existence Conjecture has had a long history in Design Theory 
since 1853 when Steiner asked about the existence of 
Steiner $(n,q,r)$ systems. Here we briefly mention a few highlights
that are relevant to our discussion here. The case $r=2$ was
proved by Wilson \cite{W1,W2,W3} in the 1970's. Around the same time,
Graver and Jurkat \cite{GJ} and Wilson \cite{W4} showed that the
divisibility condition suffices for an integral $(n,q,r,\lL)$ design,
i.e.\ an assignment of integer weights $w_Q$ to copies $Q$ of $K^r_q$
in $K^r_n$ such that $\sum \{ w_Q: e \in Q\} = \lL$ for all $e \in K^r_n$.
R\"odl \cite{R} showed the existence of approximate Steiner systems,
i.e.\ that there are edge-disjoint copies of $K^r_q$ in $K^r_n$
such that only $o(n^r)$ edges are not covered; his semi-random (nibble)
method is now an indispensable tool of modern Probabilistic Combinatorics.
Teirlinck \cite{T} was the first to show that there are {\em any} non-trivial
$(n,q,r,\lL)$ designs for arbitrary $r$. Kuperberg, Lovett and Peled \cite{KLP}
gave an alternative probabilistic proof of this result (and the existence
of many other regular combinatorial structures); their method was extended
by Lovett, Rao and Vardy \cite{LRV} to show the existence of `large sets'
of designs (for certain parameter sets). 
Glock, K\"uhn, Lo and Osthus \cite{GKLO} gave an alternative
combinatorial proof of the Existence Conjecture 
(the proof in \cite{Kexist} used a randomised algebraic construction);
they also weakened the typicality hypothesis of \cite{Kexist} (version 1)
to an extendability hypothesis, similar to that subsequently used
in \cite{Kexist} (version 2).
Furthermore, in \cite{GKLO2} they obtained analogous results
on $H$-decompositions where $H$ is any $r$-graph
and $G$ is an $r$-graph that is $H$-divisible,
i.e.\ each degree $|G(e)|$ is divisible by the gcd
of all degrees $|H(f)|$ with $|f|=|e|$.

Having discussed some hypergraph generalisations of Kirkman's result 
on triangle decompositions of $K_n$ (Steiner Triple Systems), 
let us now consider such generalisations for triangle
decompositions of $K_3(n)$ (Latin Squares).
Besides being a combinatorially natural direction,
this also has practical applications.
For example, in software testing (see \cite{H}),
a $K^r_q$-decomposition of\footnote{
For any hypergraph $H$ we write $H(n)$ for its $n$-blowup.}
$K^r_q(n)$ can be thought
of as a sequence of tests to a program taking $q$ inputs from $[n]$,
so that for every $r$ inputs all possible combinations are tested once
(so an efficient $K^r_q$-covering of $K^r_q(n)$ suffices in this context).
Another example is to a secret sharing scheme
that distributes information to $q-1$ bank clerks
so that any $r$ of them can open the safe but any $r-1$ cannot:
pick a random copy of $K^r_q$ in the decomposition, give one vertex 
to each clerk, and make the final vertex the combination for the safe.
High-dimensional permutations (also called Latin Hypercubes)
are equivalent to $K^r_{r+1}$-decompositions of $K^r_{r+1}(n)$.
In section \ref{sec:part}
we will show how the result of \cite{K2} implies an approximate
formula for the number of such decompositions, thus confirming
a conjecture of Linial and Luria \cite{LL}. The method applies in 
greater generality: as an other illustration we will give 
an approximate formula for the number of generalised Sudoku squares,
via $H$-decompositions of $H(n)$ for an auxiliary $4$-graph $H$.

In section \ref{sec:genpart} we consider a common generalisation
of the nonpartite and partite decompositions discussed above
to a generalised partite setting in which the edges of $H$ and $G$ 
have the same intersection patterns with respect to some partitions
of their vertex sets. This general setting encodes several further
problems in Design Theory. For example, Kirkman's Schoolgirl Problem
(a popular puzzle in the 19th century) asks for the construction of
a Steiner Triple System that is resolvable,
meaning that its blocks can be partitioned into perfect matchings
(sets of triples covering every vertex exactly once).
We will illustrate the generalisation to 
hypergraph decompositions given in \cite{K2}.
We will also illustrate the construction in \cite{K2} of large sets 
of designs, i.e.\ decomposition of $K^q_n$ into $(n,q,r,\lL)$ designs.
An application of the latter (see \cite{SS})
is to the following `Russian Cards' problem in information security.
From a deck of $n$ cards, we randomly deal cards so that Alice 
receives $a$ cards, Eve $e<a$ cards and Bob $b=n-a-e$ cards.
Alice wants to make a public announcement from which Bob can 
learn her cards (given the cards that he holds) 
while limiting the information that Eve receives
(e.g.\ for any card that she does not hold she should 
not learn which of Alice or Bob holds it).
A strategy for this problem can be identified with a partition of $K^a_n$,
where edges represent the possible sets of cards for Alice,
and Alice announces to which part her actual set belongs.
An optimal (minimum number of parts) strategy such that Bob
can learn Alice's hand corresponds to a 
partition of $K^a_n$ into Steiner $(n,a,a-e)$ systems;
furthermore, if $n>n_0(a,e)$ is large then it is secure against Eve,
as for any card $x$ that she does not hold,
among the blocks disjoint from her hand
in any of the Steiner systems, 
at least one contains $x$
and at least one does not.

We will explain the statement of the result of \cite{K2}
in section \ref{sec:theory}, and illustrate it with two new
applications in the subsequent two sections.
In section \ref{sec:col} we generalise the results on
hypergraph decomposition discussed above to decompositions
of hypergraphs where edges have colours which must be 
respected by the decomposition. As well as being combinatorially
natural, such generalisations encode other problems of Design Theory
(e.g.\ Whist Tournaments) and also fit within the large
literature on rainbow versions of classical combinatorial results,
which can encode seemingly unrelated questions (see e.g.\ \cite{MPS}).
In section \ref{sec:di} we give a different generalisation,
namely to decompositions of directed hypergraphs. 
This illustrates the following important feature 
of the result of \cite{K2}: it is fundamentally concerned
with sets of functions (which we call labelled edges),
so to apply it to sets of (unlabelled) edges (i.e.\ hypergraphs)
we must encode an edge by a suitable set of labelled edges.
This general setting has more applications, albeit at the expense of
considerable effort required in setting up the theory 
in section \ref{sec:theory}. However, this seems unavoidable,
as there are divisibility phenomena even for unlabelled coloured hypergraphs
that require labels to analyse (see \cite[section 1.5]{K2}).
In section \ref{sec:all} we give a common generalisation of the
previous results for convenient use in applications.
We conclude in section \ref{sec:end}
by discussing some directions for potential future research.

\section{Partite decompositions, hypermutations, Sudoku} \label{sec:part}

Over the next three sections we will gradually move from 
examples to the general setting. We start with this section
by illustrating some results on hypergraph decompositions 
and some of their applications discussed in introduction.
First we consider the nonpartite setting with
the typicality condition from \cite{Kexist},
which describes an $r$-graph where the common neighbourhood 
of small set of $(r-1)$-sets behaves roughly as one would 
expect in a random $r$-graph of the same density.

\begin{defn} \label{def:typ}
Suppose $G$ is an $r$-graph on $[n]$.
The density of $G$ is\footnote{We identify any hypergraph
with its edge-set, so $|G|$ is the number of edges.}
$d(G) = |G| \tbinom{n}{r}^{-1}$.
We say that $G$ is \emph{$(c,s)$-typical} if for any set $A$ 
of $(r-1)$-subsets of $V(G)$ with $|A| \le s$ we have 
$\bsize{\cap_{f \in A} G(f)} = (1 \pm |A|c) d(G)^{|A|} n$.
\end{defn}

The following result of \cite{GKLO2} (see also \cite[Theorem 1.5]{K2})
shows that any dense typical $r$-graph has an $H$-decomposition provided 
that it satisfies the necessary divisibility condition discussed above.
Henceforth we fix parameters
\[ h=2^{50q^3} \ \text{ and } \ \dD = 2^{-10^3 q^5}.\]

\begin{theo} \label{Hdecomp:typ}
Let $H$ be an $r$-graph on $[q]$
and $G$ be an $H$-divisible
$(c,h^q)$-typical $r$-graph on $[n]$,
where $n>n_0(q)$ is large,
$d(G) > 2n^{-\dD/h^q}$, $c < c_0 d(G)^{h^{30q}}$
and $c_0=c_0(q)$ is small.
Then $G$ has an $H$-decomposition.
\end{theo}

Next we set up some notation for stating
the partite analogue of the previous result.

\begin{defn} \label{def:Hblowup}
Let $H$ be an $r$-graph. 
We call an $r$-graph $G$ an $H$-blowup
if $V(G)$ is partitioned as $(V_x: x \in V(H))$
and each $e \in G$ is $f$-partite for some $f \in H$,
i.e.\ $f = \{x: e \cap V_x \ne \es \}$.

We write $G_f$ for the set of $f$-partite $e \in G$.
For $f \in H$ let $d_f(G) = |G_f| \prod_{x \in f} |V_x|^{-1}$.
We call $G$ a $(c,s)$-typical $H$-blowup if 
for any $s' \le s$ and distinct $e_1,\dots,e_{s'}$ where
each $e_j$ is $f_j$-partite for some $f_j \in \tbinom{V(H)}{r-1}$,
and any $x \in \cap_{j=1}^{s'} H(f_j)$ we have
$\bsize{V_x \cap \bigcap_{j=1}^{s'} G(e_j)}
= (1 \pm s'c) |V_x| \prod_{j=1}^{s'} d_{f_j \cup \{x\}}(G)$.

We say $G$ has a partite $H$-decomposition if 
it has an $H$-decomposition using copies of $H$ 
with one vertex in each part $V_x$.

We say $G$ is $H$-balanced if for every $f \sub V(H)$
and $f$-partite $e \sub V(G)$ there is some $n_e$
such that $|G_{f'}(e)|=n_e$ for all $f \sub f' \in H$.
\end{defn}

Note in particular that the $H$-balance condition for $e=f=\es$ 
implies equality of all $|G_{f'}|$ with $f' \in H$.
If $G$ has a partite $H$-decomposition then $G$ must be
$H$-balanced; the following result (\cite[Theorem 1.7]{K2})
shows the converse for typical $H$-blowups.

\begin{theo} \label{Hdecomp:partite1}
Let $H$ be an $r$-graph on $[q]$
and $G$ be an $H$-balanced $(c,h^q)$-typical
$H$-blowup on $(V_x: x \in V(H))$,
where each $n/h \le |V_x| \le n$
for some large $n > n_0(q)$ and $d_f(G) > d > 2n^{-\dD/h^q}$ 
for all $f \in H$ and $c < c_0 d^{h^{30q}}$,
where $c_0=c_0(q)$ is small.
Then $G$ has a partite $H$-decomposition.
\end{theo}

In the previous result, we can not only show
that $G$ has a partite $H$-decomposition,
but also give an approximate formula for the
number of such decompositions. We will show some 
applications of this when $G$ is a complete $H$-blowup.
We start by considering the upper bound,
which comes from the following result of Luria \cite{L}.

\begin{theo} \label{luria}
Let $R$ be fixed and $D=D(N) \to \infty$ as $N \to \infty$.
Suppose $A$ is an $R$-graph on $N$ vertices
such that all vertex degrees are\footnote{
The statement in \cite{L} has $D$ here,
but the proof works with $D+o(D)$.}  
$D+o(D)$ and all pair degrees are $o(D)$.
Then the number of perfect matchings in $A$ 
is at most $(De^{1-R} + o(D))^{N/R}$.
\end{theo}

When applying Theorem \ref{luria} to the setting of 
Theorem \ref{Hdecomp:partite1}, we consider the auxiliary $R$-graph 
$A$ on $V(A)=E(G)$ where edges correspond to copies of $H$, 
so $N=|G|$ and $R=|H|$. If we let $G=H(n)$ be the complete $H$-blowup
of size $n$ then $N=|H|n^r$ and the degree conditions of Theorem \ref{luria}
hold with $D = n^{q-r}$. In fact, all pair degrees are at most $n^{q-r-1}$.
We deduce that the number of $H$-decompositions 
of $H(n)$ is at most $( (e^{1-|H|}+o(1)) n^{q-r} )^{n^r}$.
We will show below how a matching lower bound follows
from Theorem \ref{Hdecomp:partite1}. Before doing so,
we discuss two applications.

First we consider the number $N_r(n)$
of $r$-dimensional permutations of order $n$,
which is also the number of $K^r_{r+1}$-decompositions of $K^r_{r+1}(n)$.
For $r=2$ (Latin squares), Van Lint and Wilson \cite[Theorem 17.3]{vLW}
obtained the approximate formula $N_2(n) = (n/e^2 + o(n))^{n^2}$; 
this was a short deduction from two celebrated breakthroughs on permanents
(the proof of the Van der Waerden Conjecture by
Falikman and by Egorychev and of the Minc Conjecture by Bregman).
The upper bound can be obtained more simply by entropy inequalities,
by which means Linial and Luria \cite{L} showed 
$N_r(n) \le (n/e^r + o(n))^{n^r}$, and Luria obtained
the more general result in Theorem \ref{luria}.
However, the lower bound argument appeared not to generalise,
even from Latin squares to Steiner Triple Systems,
for which the approximate formula was a conjecture of 
Wilson \cite{W6}, proved in \cite{Kcount}.
In \cite{K2} we established the lower bound,
thus giving the following approximate formula.

\begin{theo}
The number of $r$-dimensional permutations 
of order $n$ is $(n/e^r + o(n))^{n^r}$.
\end{theo}

Our second application is to the number 
of generalised Sudoku squares, which are
Latin squares of order $n^2$ partitioned into
$n$ by $n$ subsquares each of which uses every symbol once
(the usual Sudoku squares have $n=3$).
We encode these by the $4$-graph $H$ with
$V(H) = \{x_1,x_2,y_1,y_2,z_1,z_2\}$ and
$E(H) = \{ x_1x_2y_1y_2, x_1x_2z_1z_2, 
y_1y_2z_1z_2, x_1y_1z_1z_2 \}$.
Then an $H$-decomposition of the complete $n$-blowup of $H$
can be viewed as a Sudoku square, 
where we represent rows by pairs $(a_1,a_2)$,
columns by $(b_1,b_2)$, symbols by $(c_1,c_2)$
and boxes by $(a_1,b_1)$; a copy of $H$
with vertices $\{a_1,a_2,b_1,b_2,c_1,c_2\}$
represents a cell in row $(a_1,a_2)$
and column $(b_1,b_2)$ with symbol $(c_1,c_2)$.
The following estimate then follows from the estimate
for general $H$ given below.

\begin{theo}
The number of Sudoku squares with $n^2$ boxes 
of order $n$ is $(n^2/e^3 + o(n^2))^{n^4}$.
\end{theo}

We conclude this section with the general formula
that implies the two examples discussed above.

\begin{theo}
For any $r$-graph $H$ on $[q]$,
the number of $H$-decompositions of $H(n)$
is $( (e^{1-|H|}+o(1)) n^{q-r} )^{n^r}$.
\end{theo}

\nib{Proof.}
The upper bound comes from Theorem \ref{luria}
applied to the auxiliary $R$-graph $A$ described above
(following the statement of Theorem \ref{luria}).
For the lower bound, we consider the random greedy matching process,
in which we construct a sequence of vertex-disjoint edges 
$e_0,e_1,\dots$ in $A$ and subgraphs $A_0,A_1,\dots$,
where $A_0=A$, $e_i$ is a uniformly random edge of $A_i$,
and $A_{i+1}$ is obtained from $A_i$ by deleting the vertices
of $e_i$ and all edges that intersect $e_i$.
We will estimate the number of runnings of this process,
stopped at some subgraph $A_t$ which is quite sparse,
but sufficiently dense and typical that Theorem \ref{Hdecomp:partite1}
applies to show that $A_t$ has a perfect matching.
This will give a lower bound on the number of perfect matchings of $A$,
i.e.\ $H$-decompositions of $H(n)$, which matches Luria's upper bound.

Bennett and Bohman \cite{BB} showed if $A$ 
is a $D$-regular $R$-graph on $N$ vertices with 
all pair degrees at most $L = o(D\log^{-5} N)$
then whp\footnote{
We say that an event $E$ holds \emph{with high probability} (whp) 
if $\mb{P}(E) = 1-e^{-\Omega(n^c)}$ for some $c>0$ as $n \to \infty$; 
by union bounds we can assume that any specified 
polynomial number of such events all occur. }
the process persists until the proportion of
uncovered vertices is at most $(L/D)^{1/2(R-1) + o(1)}$.
(Their proof applies verbatim under the weaker
assumption that all vertex degrees are $D \pm \sqrt{DL}$.)
Here we have $L/D = n^{-1}$ and $R=|H|$, 
so we could run the process until the uncovered
proportion is e.g.\ $n^{-1/2|H|}$, but we stop it 
when the remaining $r$-graph $G_t=V(A_t)$ has density $d = 3n^{-\dD/h^q}$.
Furthermore, one can show that whp throughout the process
the $r$-graphs $G_i=V(A_i)$ are $(c,h^q)$-typical $H$-blowups 
with $c < c_0 d^{h^{30q}}$ (similar lemmas in the nonpartite setting
are well-known, see e.g.\ \cite{BFL}).
Then Theorem \ref{Hdecomp:partite1} can be applied to $G_t$,
and we have a good estimate for the number of choices
at each step of the process: at step $i$
when all densities $d_f(G_i)$ with $f \in H$
are $d(i) = 1 - in^{-r}$
there are $(1 \pm 2|H|c) d(i)^{|H|} n^q$ 
edges of $A_i$ (i.e.\ copies of $H$ in $G_i$).

Given the above results, a simple counting argument now gives the 
required lower bound on the number of $H$-decompositions of $H(n)$.
For $0 \le j \le j' \le t$, let us say that a running of the process 
from $A_0,\dots,A_{j'}$ is $j$-good if $G_i$ is $(c,h^q)$-typical 
for $1 \le i \le j$. Let $R^j_{j'}$ be the number of such runnings.
Then $R^j_{j+1} / R^j_j = (1 \pm 2|H|c) d(j)^{|H|} n^q$ by typicality
and $R^{j+1}_{j+1} / R^j_{j+1} = 1 \pm c$ (say) as whp
typicality does not first fail at step $j+1$.
Multiplying these estimates, the number of $t$-good runnings
is $R^t_t = \prod_{j=0}^t ( (1 \pm 3|H|c) d(j)^{|H|} n^q )$.
By Theorem \ref{Hdecomp:partite1}, each $t$-good 
running can be completed to an $H$-decomposition of $H(n)$.
We obtain a lower bound on the number of $H$-decompositions of $H(n)$
by dividing $R^t_t$ by an upper bound of
$\prod_{j=0}^t (n^r-j) = \prod_{j=0}^t (d(j) n^r)$
on  the number of runnings giving rise to any fixed decomposition.
A short calculation using Stirling's estimate on factorials
gives the claimed lower bound
$\prod_{j=0}^t ( (1 \pm 3|H|c) d(j)^{|H|-1} n^{q-r} )
= ( (e^{1-|H|}+o(1)) n^{q-r} )^{n^r}$. \qed

\section{Generalised partite decompositions} \label{sec:genpart}

In this section we state and give applications of 
a result that generalises both the nonpartite and partite
decomposition results of the previous section
to the generalised partite setting of the definition below
(which is followed by some explanatory remarks).

\begin{defn} \label{def:HPblowup}
Let $H$ be an $r$-graph on $[q]$ and 
$\mc{P}=(P_1,\dots,P_t)$ be a partition of $[q]$.
Let $G$ be an $r$-graph and
$\mc{P}'=(P'_1,\dots,P'_t)$ be a partition of $V(G)$.
We say $G$ has a $\mc{P}$-partite $H$-decomposition if 
it has an $H$-decomposition using copies $\phi(H)$ 
of $H$ with all $\phi(P_i) \sub P'_i$.

For $S \sub [q]$ the $\mc{P}$-index of $S$ is 
$i_{\mc{P}}(S) = (|S \cap P_1|,\dots,|S \cap P_t|)$;
similarly, we define the $\mc{P}'$-index of subsets
of $V(G)$, and also refer to both as the `index'.

For $\iB \in \mb{N}^t$ we let $H_\iB$ and $G_\iB$ be the
edges in $H$ and $G$ of index $\iB$.
Let $I=I(H) = \{\iB: H_\iB \ne \es \}$.
We call $G$ an $(H,\mc{P})$-blowup 
if $G_\iB \ne \es \Ra i \in I$. 

For $e \sub V(G)$ we define 
the degree vector $G_I(e) \in \mb{N}^I$
by $G_I(e)_\iB=|G_\iB(e)|$ for $\iB \in I$.
Similarly, for $f \sub [q]$ we define
$H_I(f)$ by $H_I(f)_\iB=|H_\iB(f)|$.
For $\iB' \in \mb{N}^t$ let $H^I_{\iB'}$
be the subgroup of $\mb{Z}^I$ generated 
by $\{ H_I(f): i_{\mc{P}}(f) = \iB' \}$.
We say $G$ is $(H,\mc{P})$-divisible 
if $G_I(e) \in H^I_{\iB'}$ whenever $i_{\mc{P}'}(e)=\iB'$.

For $\iB \in \mb{N}^t$ let 
$d_\iB(G) = |G_\iB| \prod_{j \in [t]} \tbinom{|P'_j|}{i_j}^{-1}$.
We call $G$ a $(c,s)$-typical $(H,\mc{P})$-blowup 
if for any $s' \le s$, 
$\{ f_1,\dots,f_{s'} \} \sub \tbinom{V(G)}{r-1}$, 
$j \in [t]$ we have\footnote
{Let $\{e_1,\dots,e_t\}$ be the standard basis of $\mb{Z}^t$.} 
$\bsize{P'_j \cap \bigcap_{k=1}^{s'} G(f_k)}
= (1 \pm s'c) |P'_j| \prod_{k=1}^{s'} d_{i(f_k)+e_j}(G)$.
\end{defn}

The simplest examples of the previous definition
are given by the trivial partitions 
with $t=1$ (non-partite decompositions)
or $t=q$ (partite decompositions).
The latter is instructive for understanding the
divisibility condition. We will illustrate it
in the case that $H$ is a (graph) triangle on $[3]$,
with parts $P_i = \{i\}$ for $i \in [3]$
and $G$ is a tripartite graph
with parts $P'_i$ for $i \in [3]$.
Then $I = \{ i^1, i^2, i^3 \}$ with 
$\iB^1 = (1,1,0)$, $\iB^2 = (1,0,1)$, $\iB^3 = (0,1,1)$. 
For each $\iB \in I$ we have $G(\es)_\iB = |G_\iB|$
and $H(\es)_\iB = |H_\iB|=1$, so the 0-divisibility condition 
is that the three bipartite pieces of $G$ 
all have the same number of edges.
For the 1-divisibility condition, we note that
$H(1)_{\iB^1}=H(1)_{\iB^2}=1$, $H(1)_{\iB^3}=0$ and
$G(x_1)_\iB = |G_\iB(x_1)|$ for $x_1 \in P'_1$,
so we require every vertex in $P'_1$ to have
equal degrees into $P'_2$ and $P'_3$
(and similarly for each part).
The 2-divisibility condition is trivially satisfied,
so this completes the description.
Our final remark on Definition \ref{def:HPblowup} is that
the typicality condition is a direct generalisation 
of that in Definition  \ref{def:Hblowup}, allowing 
the possibility that both sides are zero
if some $i(f_k)+e_j \notin I$.

Next we state a decomposition result
in the generalised partite setting
(a case of \cite[Theorem 7.8]{K2});
the case $\mc{P}=([q])$
implies Theorem \ref{Hdecomp:typ} 
and the case $\mc{P}=(\{1\},\dots,\{q\})$ 
implies Theorem \ref{Hdecomp:partite1}.

\begin{theo} \label{HPdecomp:typ}
Let $H$ be an $r$-graph on $[q]$ and 
$\mc{P}=(P_1,\dots,P_t)$ be a partition of $[q]$.
Let $n > n_0(q)$, $d > 2n^{-\dD/h^q}$ 
and $c < c_0 d^{h^{30q}}$, where $c_0=c_0(q)$ is small.
Suppose $G$ is an $(H,\mc{P})$-divisible
$(c,h)$-typical $(H,\mc{P})$-blowup 
wrt $\mc{P}'=(P'_1,\dots,P'_t)$,
such that each $n/h \le |P'_i| \le n$
and $d_i(G)>d$ for all $i \in I(H)$.
Then $G$ has a $\mc{P}$-partite $H$-decomposition.
\end{theo}

In the remainder of this section we give
two applications of the following simplified
version of the preceding result (the case that $G$ is complete).

\begin{theo} \label{HPdecomp:complete}
Let $H$ be an $r$-graph on $[q]$ and 
$\mc{P}=(P_1,\dots,P_t)$ be a partition of $[q]$.
Suppose $G$ is an $(H,\mc{P})$-divisible
complete $(H,\mc{P})$-blowup wrt $\mc{P}'=(P'_1,\dots,P'_t)$
such that each $n/h \le |P'_i| \le n$ with $n > n_0(q)$.
Then $G$ has a $\mc{P}$-partite $H$-decomposition.
\end{theo}

As our first application we reprove the result of \cite{RW}
in the case that $n$ is large on the existence of
resolvable Steiner Triple Systems (for a hypergraph
generalisation see \cite[Theorem 7.9]{K2}).

\begin{theo}
Suppose $n = 6k+3$ with $k \in \mb{N}$ is large.
Then there is a resolvable Steiner Triple System of order $n$.
\end{theo}

\nib{Proof.}
Let $H=K_4$ be the complete graph on $4$ vertices,
with $V(H)=[4]$ partitioned as $\mc{P}=(P_1,P_2)$,
where $P_1=[3]$ and $P_2=\{4\}$. 
Let $P'_1$ and $P'_2$ be disjoint sets 
with $|P'_1|=n$ and $|P'_2|=(n-1)/2$.
Let $G$ be the graph with $V(G)=P'_1 \cup P'_2$ whose edges
are all pairs in $P'_1 \cup P'_2$ not contained in $P'_2$.
Then $G$ is a complete $(H,\mc{P})$-blowup. 
 
We claim that a resolvable Steiner Triple System of order $n$
is equivalent to a $\mc{P}$-partite $H$-decomposition of $G$.
To see this, suppose first that we have some 
$\mc{P}$-partite $H$-decomposition $\mc{H}$ of $G$.
This means that $\mc{H}$ partitions $E(G)$,
and each $\phi(H) \in \mc{H}$ has
$\phi([3]) \sub P'_1$ and $\phi(4) \in P'_2$.
Then $\mc{T} := \{ \phi(H-4): \phi(H) \in \mc{H} \}$
is a triangle decomposition of the complete graph on $P_1$,
i.e.\ a Steiner Triple System of order $n$.
We can partition $\mc{H}$ as $(\mc{H}_y: y \in P'_2)$,
where each $\mc{H}_y = \{ \phi(H): \phi(4)=y \}$.
Note that each $T_y = \{ \phi([3]): \phi(H) \in \mc{H}_y \}$
is a perfect matching on $P_1$; indeed, for each $x \in P_1$,
as $\mc{H}$ partitions $E(G)$, there is a unique 
$\phi(H) \in \mc{H}$ containing $xy$, and then
$\phi([3])$ is the unique triple in $T_y$ containing $x$.
Thus $\mc{T}$ is a resolvable Steiner Triple System.
Conversely, the same construction shows that
any resolvable Steiner Triple System gives
rise to a $\mc{P}$-partite $H$-decomposition of $G$.
Indeed, given a Steiner Triple System $\mc{T}$ on $P_1$
partitioned into perfect matchings, we arbitrarily
label the perfect matchings as $(T_y: y \in P'_2)$
and form a $\mc{P}$-partite $H$-decomposition of $G$
by taking all $\phi(H)$ with 
$\phi([3]) \in T_y$ and $\phi(4)=y$ for some $y \in P'_2$.
This proves the claim.

To complete the proof of the theorem, we show that 
Theorem \ref{HPdecomp:complete} applies to give
a $\mc{P}$-partite $H$-decomposition of $G$.
In the notation of Definition \ref{def:HPblowup},
we have $I = I(H) = \{ (2,0), (1,1) \}$ 
and need to show that
$G_I(e) \in H^I_{\iB'}$ whenever $i_{\mc{P}'}(e)=\iB'$.
First we consider $i_{\mc{P}'}(e)=(0,0)$, i.e.\ $e=\es$.
We have $H_I(\es) = (3,3)$, as $H$ contains $3$ edges of
each of the indices $(2,0)$ and $(1,1)$.
Thus $H^I_{(0,0)} \le \mb{Z}^2$ is generated by $(3,3)$.
We have $G_I(\es) = ( \tbinom{n}{2}, \tbinom{n}{2} )$,
as $G$ contains $\tbinom{n}{2}$ edges inside $P'_1$
and $\tbinom{n}{2}$ edges between $P'_1$ and $P'_2$.
As $3 \mid n$ we have $G_I(\es) \in H^I_{(0,0)}$.

Next we consider $i_{\mc{P}'}(e)=(1,0)$, i.e.\ $e \in P'_1$.
We have $i_{\mc{P}}(f)=(1,0)$ iff $f \in [3]$,
and for any such $f$ we have $H_I(f) = (2,1)$,
as $f$ is contained in $2$ edges of index $(2,0)$
and $1$ edge of index $(1,1)$.
Thus $H^I_{(1,0)} \le \mb{Z}^2$ is generated by $(2,1)$.
We have $G_I(e) = ( n-1, (n-1)/2 )$,
as $e$ has degree $n-1$ in $P'_1$
and degree $(n-1)/2$ in $P'_2$.
As $n$ is odd, $G_I(e) \in H^I_{(1,0)}$.
The only remaining non-trivial case is that
$i_{\mc{P}'}(e)=(0,1)$, i.e.\ $e \in P'_2$.
We have $i_{\mc{P}}(f)=(0,1)$ iff $f=4$,
and $H_I(4) = (0,3)$, as $f$ is contained in no edges 
of index $(2,0)$ and $3$ edges of index $(1,1)$.
Thus $H^I_{(0,1)} \le \mb{Z}^2$ is generated by $(0,3)$.
We have $G_I(e) = ( 0, n)$, as $e$ has degree $0$ in $P'_2$
and degree $n$ in $P'_1$. 
As $3 \mid n$ we have $G_I(e) \in H^I_{(1,0)}$. \qed

\medskip

Our second application is to reprove the existence
of large sets of Steiner Triple Systems for large $n$
(due to Lu, completed by Teirlinck, see \cite{T2});
see \cite[Theorem 1.2]{K2} for the hypergraph version.

\begin{theo}
Suppose $n$ is large and $1$ or $3$ mod $6$.
Then $K^3_n$ can be decomposed into Steiner Triple Systems.
\end{theo}

\nib{Proof.}
Let $H=K_4$ be the complete $3$-graph on $4$ vertices,
with $V(H)=[4]$ partitioned as $\mc{P}=(P_1,P_2)$,
where $P_1=[3]$ and $P_2=\{4\}$. 
Let $P'_1$ and $P'_2$ be disjoint sets 
with $|P'_1|=n$ and $|P'_2|=n-2$.
Let $G$ be the $3$-graph with $V(G)=P'_1 \cup P'_2$ 
whose edges are all triples 
$e \sub P'_1 \cup P'_2$ with $|e \cap P'_1| \ge 2$.
Then $G$ is a complete $(H,\mc{P})$-blowup.

We claim that a decomposition of $K^3_n$ into Steiner Triple Systems
is equivalent to a $\mc{P}$-partite $H$-decomposition of $G$.
To see this, suppose we have some 
$\mc{P}$-partite $H$-decomposition $\mc{H}$ of $G$.
We can partition $\mc{H}$ as $(\mc{H}_y: y \in P'_2)$,
where each $\mc{H}_y = \{ \phi(H): \phi(4)=y \}$.
Note that each $T_y = \{ \phi([3]): \phi(H) \in \mc{H}_y \}$
is a Steiner Triple System on $P_1$; indeed, for each pair $xx'$ in 
$P_1$, as $\mc{H}$ partitions $E(G)$, there is a unique 
$\phi(H) \in \mc{H}$ containing $xx'y$, and then
$\phi([3])$ is the unique triple in $T_y$ containing $xx'$.
Furthermore, each triple in $P'_1$ belongs to exactly one
element of $\mc{H}$, and so to exactly one $T_y$.
Thus $\{T_y: y \in P'_2\}$ is a decomposition 
of $K^3_n$ into Steiner Triple Systems.
Conversely, the same construction converts 
any decomposition of $K^3_n$ into Steiner Triple Systems
into a $\mc{P}$-partite $H$-decomposition of $G$.

To complete the proof of the theorem, we show that 
Theorem \ref{HPdecomp:complete} applies to give
a $\mc{P}$-partite $H$-decomposition of $G$.
We have $I = I(H) = \{ (3,0), (2,1) \}$
and need to show that
$G_I(e) \in H^I_{\iB'}$ whenever $i_{\mc{P}'}(e)=\iB'$.
First we consider $i_{\mc{P}'}(e)=(a,0)$ with $0 \le a \le 2$.
For any $f \sub V(H)$ with $i_{\mc{P}}(f)=(a,0)$ we have
$H_I(f) = (1,3-a)$, as $f$ is contained in $1$ edge of $H$
with index $(3,0)$ and $3-a$ edges of $H$ with index $(2,1)$.
Thus $H^I_{(a,0)} \le \mb{Z}^2$ is generated by $(1,3-a)$.
We have $G_I(e) = ( \tbinom{n-a}{3-a}, (3-a) \tbinom{n-a}{3-a} )$,
as $e$ is contained in $\tbinom{n-a}{3-a}$
edges of $G$ with index $(3,0)$
and $\tbinom{n-a}{2-a} (n-2) = (3-a) \tbinom{n-a}{3-a}$
edges of $G$ with index $(2,1)$.
Therefore $G_I(e) \in H^I_{(a,0)}$.

Next consider $i_{\mc{P}'}(e)=(0,1)$, i.e.\ $e \in P'_2$.
We have $i_{\mc{P}}(f)=(0,1)$ iff $f=4$,
and $H_I(4) = (0,3)$, as $4$ is contained in $0$ edges
of index $(3,0)$ and $3$ edges of index $(2,1)$.
Thus $H^I_{(0,1)} \le \mb{Z}^2$ is generated by $(0,3)$.
We have $G_I(e) = ( 0, \tbinom{n}{2} )$,
as $e$ is contained in no edges of $G$ with index $(3,0)$
and $\tbinom{n}{2}$ edges of $G$ with index $(2,1)$.
As $3 \mid \tbinom{n}{2}$ we have $G_I(e) \in H^I_{(0,1)}$.

The only remaining non-trivial case is $i_{\mc{P}'}(e)=(1,1)$.
We have $i_{\mc{P}}(f)=(1,1)$ iff $f=a4$ for some $a \in [3]$.
Then $H_I(f) = (0,2)$, as $f$ is contained in $0$ edges
of index $(3,0)$ and $2$ edges of index $(2,1)$.
Thus $H^I_{(1,1)} \le \mb{Z}^2$ is generated by $(0,2)$.
We have $G_I(e) = ( 0, n-1 )$,
as $e$ is contained in no edges of $G$ with index $(3,0)$
and $n-1$ edges of $G$ with index $(2,1)$.
As $n$ is odd, $G_I(e) \in H^I_{(1,1)}$. \qed

\section{General theory} \label{sec:theory}

In this section we state the main result of \cite{K2},
from which all the other results in this paper follow.
Most of the section will be occupied with preparatory 
definitions for the statement of the result,
which we will illustrate with the following running example.
Consider a graph $G$ with $V(G)=[n]$ partitioned as $(V_1,V_2)$,
where there are no edges within $V_2$, edges within $V_1$ are red, 
and edges between $V_1$ and $V_2$ are blue or green.
When does $G$ have a decomposition into rainbow triangles?

\subsection{Labelled complexes and embeddings}

All decomposition problems that fit in our general framework
are encoded by labelled complexes, which are sets of functions
(which we think of as labelled edges) closed under taking restriction; 
this is analogous to (simplicial) complexes, which are sets of sets
closed under taking subsets. 

\begin{defn} \label{def:complex} \

We call $\Phi=(\Phi_B: B \sub R)$ an $R$-system on $V$ 
if $\phi:B \to V$ is injective for each $\phi \in \Phi_B$.

We call $\Phi$ an $R$-complex
if whenever $\phi \in \Phi_B$ and $B' \sub B$ 
we have $\phi\mid_{B'} \in \Phi_{B'}$.

Let $\Phi^\circ_B = \{ \phi(B): \phi \in \Phi_B \}$,
$\Phi^\circ_j = \bigcup \{ \Phi^\circ_B: B \in \tbinom{R}{j} \}$,
$\Phi_j = \bigcup \{ \Phi_B: B \in \tbinom{R}{j} \}$,
$V(\Phi) = \Phi^\circ_1$ and
$\Phi^\circ = \bigcup \{ \Phi^\circ_B: B \sub R\}$.
\end{defn}

To apply Definition \ref{def:complex}
in our example we take $V = V(G)$, $R = [3]$ and for each
$B \sub [3]$ we let $\Phi_B$ consist of all injections
$\phi:B \to V$ with $\phi(B \cap \{1,2\}) \sub V_1$
and $\phi(B \cap \{3\}) \sub V_2$: we also call $\Phi$
the complete $(\{1,2\},3)$-partite $[3]$-complex wrt $(V_1,V_2)$.
We think of $\phi \in \Phi_3$ as an embedding of the triangle on $[3]$
where $12$ is red, $13$ is blue and $23$ is green.
It is useful to consider all such embeddings,
even though the only ones that can appear in a
decomposition of $G$ are those that are contained in $G$ 
with $\phi(12)$ red, $\phi(13)$ blue and $\phi(23)$ green.

Next we consider the functional analogue 
of the subgraph notion for hypergraphs.
Just as an embedding of a hypergraph $H$ 
in a hypergraph $G$ is an injection
from $V(H)$ to $V(G)$ taking edges to edges,
an embedding of labelled complexes is an injection 
taking labelled edges to labelled edges.

\begin{defn}
Let $H$ and $\Phi$ be $R$-complexes.
Suppose $\phi:V(H) \to V(\Phi)$ is injective.
We call $\phi$ a $\Phi$-embedding of $H$
if $\phi \circ \psi \in \Phi$ for all $\psi \in H$.
\end{defn}

In our example, $\Phi$ is as above, and $H$ is the complete 
$(\{1,2\},3)$-partite $[3]$-complex wrt $(\{1,2\},3)$,
i.e.\ each $H_B$ with $B \sub [3]$ consists of all injections
$\phi:B \to [3]$ with $\phi(B \cap \{1,2\}) \sub \{1,2\}$
and $\phi(B \cap \{3\}) \sub \{3\}$.
We think of an edge $e$ of the triangle on $[3]$
as being encoded by the set of labelled edges of $H$
with image $e$, thus $12$ is encoded
by $\{(1 \mapsto 1, 2 \mapsto 2),
(1 \mapsto 2, 2 \mapsto 1)\}$,
$13$ by $\{(1 \mapsto 1, 3 \mapsto 3),
(2 \mapsto 1, 3 \mapsto 3)\}$,
and $23$ by $\{(2 \mapsto 2, 3 \mapsto 3),
(1 \mapsto 2, 3 \mapsto 3)\}$.
If $\phi$ is a $\Phi$-embedding of $H$
we encode the edges of the triangle on $\phi([3])$
by the corresponding sets of labelled edges:
$12$ by $\{(1 \mapsto \phi(1), 2 \mapsto \phi(2)),
(1 \mapsto \phi(2), 2 \mapsto \phi(1))\}$,
$13$ by $\{(1 \mapsto \phi(1), 3 \mapsto \phi(3)),
(2 \mapsto \phi(1), 3 \mapsto \phi(3))\}$,
and $23$ by $\{(2 \mapsto \phi(2), 3 \mapsto \phi(3)),
(1 \mapsto \phi(2), 3 \mapsto \phi(3))\}$.

\subsection{Extensions and extendability} \label{sec:ext}

Next we will formulate our extendability condition.

\begin{defn}
Let $R(S)$ be the $R$-complex of all 
partite maps from $R$ to $R \times S$, i.e.\ 
whenever $i \in B \sub R$ and $\psi \in R(S)_B$ 
we have $\psi(i)=(i,x)$ for some $x \in S$.
If $S=[s]$ we write $R(S)=R(s)$.
\end{defn}

\begin{defn} \label{def:ext}
Suppose $J \sub R(S)$ is an $R$-complex
and $F \sub V(J)$. Define $J[F] \sub R(S)$ by
$J[F] = \{ \psi \in J: Im(\psi) \sub F \}$.
Suppose $\phi$ is a $\Phi$-embedding of $J[F]$.
We call $E=(J,F,\phi)$ a $\Phi$-extension of rank $s=|S|$.
We write $X_E(\Phi)$ for the set or number of 
$\Phi$-embeddings of $J$ that restrict to $\phi$ on $F$.
We say $E$ is $\oO$-dense (in $\Phi$)
if $X_E(\Phi) \ge \oO |V(\Phi)|^{v_E}$,
where $v_E := |V(J) \sm F|$.
We say $\Phi$ is $(\oO,s)$-extendable
if all $\Phi$-extensions of rank $s$ are $\oO$-dense.
\end{defn}

In our example, we could consider extending some fixed 
rainbow triangle to an octahedron in which every 
triangle is rainbow. To implement this in the preceding two
definitions, we let $J = [3](2)$ and $F = [3] \times \{1\}$.
We identify $F$ with $[3]$ by identifying each $(i,1)$ with $i$.
Then $J[F]_B = \{id_B\}$ for $B \sub [3]$ and 
$\phi$ is a $\Phi$-embedding of $J[F]$ iff $\phi \in \Phi_3$.
We think of $Im(\phi)$ as our fixed rainbow triangle,
which has $2$ vertices in $V_1$ and $1$ vertex in $V_2$.
Now consider any $\phi^+ \in X_E(\Phi)$ where $E=(J,F,\phi)$,
i.e.\ $\phi^+$ is a $\Phi$-embedding of $J$ 
that restricts to $\phi$ on $F$.
For each $i \in [3]$ we have $(i \mapsto (i,2)) \in J_1$,
so $(i \mapsto \phi^+((i,2)) \in \Phi_1$;
thus $\phi^+((i,2)) \in V_1$ if $i \in [2]$
or $\phi^+((i,2)) \in V_2$ if $i=3$.
We think of $\{\phi^+((i,1)), \phi^+((i,2))\}$ for $i \in [3]$
as the opposite vertices of an octahedron extending
the fixed triangle $Im(\phi)$.
(We do not yet consider the colours; these will
come into play when we consider Definition \ref{def:ext+}.)
We have $X_E(\Phi) = (|V_1|-2)(|V_1|-3)(|V_2|-1)$,
as we choose $2$ new vertices in $V_1$ and $1$ in $V_2$,
so $E$ is $\OO(1)$-dense if $|V_1|$ and $|V_2|$ are both $\OO(n)$.

Next we augment our extendability condition
to allow for various restrictions 
(coloured edges in our example).

\begin{defn} \label{def:ext+}
Let $\Phi$ be an $R$-complex and
$\Phi' = (\Phi^t: t \in T)$ 
with each $\Phi^t \sub \Phi$.
Let $E=(J,F,\phi)$ be a $\Phi$-extension 
and $J' = (J^t: t \in T)$ for some
mutually disjoint $J^t \sub J \sm J[F]$;
we call $(E,J')$ a $(\Phi,\Phi')$-extension.
If $|T|=1$ we identify $\Phi' \sub \Phi$ with $(\Phi')$.

We write $X_{E,J'}(\Phi,\Phi')$ for the set or number 
of $\phi^+ \in X_E(\Phi)$ with $\phi^+ \circ \psi \in \Phi^t_B$ 
whenever $\psi \in J^t_B$ and $\Phi^t_B$ is defined.
We say $(E,J')$ is $\oO$-dense in $(\Phi,\Phi')$
if $X_{E,J'}(\Phi,\Phi') \ge \oO |V(\Phi)|^{v_E}$.

We say $(\Phi,\Phi')$ is $(\oO,s)$-extendable 
if all $(\Phi,\Phi')$-extensions of rank $s$ 
are $\oO$-dense in $(\Phi,\Phi')$.

For $G'=(G^t: t \in T)$ with each $G^t \sub \Phi^\circ$
and $J'$ as above we write 
$X_{E,J'}(\Phi,G)=X_{E,J'}(\Phi,\Phi')$, where
$\Phi' = (\Phi^t: t \in T)$ with each
$\Phi^t = \{\phi \in \Phi: Im(\phi) \in G^t\}$.

We say that $(\Phi,G')$ is $(\oO,s)$-extendable
if $(\Phi,\Phi')$ is $(\oO,s)$-extendable.
\end{defn}

We continue the above example of extending a fixed
rainbow triangle to an octahedron of rainbow triangles.
We continue to ignore colours and first consider how
the preceding definition can ensure that 
the octahedron is a subgraph of $G$.
Indeed, if $\phi^+ \in X_{E,J \sm J[F]}(\Phi,\Phi')$ with 
$\Phi'=\{\phi \in \Phi: Im(\phi) \in G\}$
then $\Phi'_B$ is only defined when $|B|=2$,
and for all $\psi \in J_2 \sm J[F]$
we have $\phi^+ \circ \psi \in \Phi'$,
i.e.\ $\phi^+(Im(\psi)) \in G$, as required.
We also note for future reference that if for some $r$
we have all $\Phi^t \sub \Phi_r$ then when checking
extendability we can assume $J' \sub J_r \sm J[F]$.

To implement colours, we let $T = \{12,13,23\}$, and for $t \in T$
let $G^t$ be the set of edges of $G$ of the appropriate colour
(red if $t=12$, blue if $t=13$, green if $t=23$),
$\Phi^t = \{\phi \in \Phi: Im(\phi) \in G^t\}$
and $J^t = J_t \sm J[F]$ for $t \in T$.
If $\phi^+ \in X_{E,J'}(\Phi,\Phi')$ then
for each $t \in T$, $\psi \in J_t \sm J[F]$
we have $\phi^+(Im(\psi)) \in G^t$, as required.
The extendability condition says that there are at least $\oO n^3$ 
such octahedra of rainbow triangles containing $\phi$ 
(and similarly for any other extension of bounded size).

\subsection{Adapted complexes}

A common feature of the decomposition results obtained from 
our main theorem is that they are implemented by a labelled complex
equipped with a permutation group action, and the decomposition 
respects the orbits of the action, as in the following definitions.

\begin{defn} \label{def:perm} 
Suppose $\Ss$ is a permutation group on $R$.
For $B,B' \sub R$ we write 
$\Ss^{B'}_B = \{ \sS\mid_B: \sS \in \Ss, \sS(B)=B' \}$,
$\Ss^{B'} = \cup_B \Ss^{B'}_B$ and
$\Ss^\le = \cup_{B,B'} \Ss^{B'}_B$.
\end{defn}

\begin{defn} \label{def:adapt} 
Suppose $\Phi$ is an $R$-complex
and $\Ss$ is a permutation group on $R$.
For $\sS \in \Ss$ and $\phi \in \Phi_{\sS(B)}$
let $\phi \sS = \phi \circ \sS\mid_B$.
We say $\Phi$ is $\Ss$-adapted if 
$\phi \sS \in \Phi$ for any $\phi \in \Phi$, $\sS \in \Ss$.
\end{defn}

\begin{defn} \label{def:orbit} 
For $\psi \in \Phi_B$ with $B \sub R$ 
we define the orbit of $\psi$ by
$\psi\Ss := \psi\Ss^B = \{ \psi\sS: \sS \in \Ss^B \}$.
We denote the set of orbits by $\Phi/\Ss$.
We write $Im(O)=Im(\psi)$ for $\psi \in O \in \Phi/\Ss$.
\end{defn}

\begin{defn} 
Let $\GG$ be an abelian group.
For $J \in \GG^{\Phi_r}$ and $O \in \Phi_r/\Ss$
we define $J^O$ by $J^O_\psi = J_\psi 1_{\psi \in O}$.
The orbit decomposition of $J$ is
$J = \sum_{O \in \Phi_r/\Ss} J^O$.
\end{defn}

The simplest example is when 
the permutation group is the entire symmetric group,
e.g.\ if $R=[3]$ and $\Ss=S_3$ then any $\phi \in \Phi_3$ 
has an orbit consisting of all six bijections from $[3]$
to $e = Im(\phi)$, which we would think of as encoding 
the edge $e$ in a $3$-graph. In our running example,
we have $\Ss = \{ id, (12) \} \le S_3$. We recall that
if $\phi$ is a $\Phi$-embedding of $H$
then the edge $\phi(12)$ of $\Phi^\circ_2$
is encoded by the labelled edges
$(1 \mapsto \phi(1), 2 \mapsto \phi(2))$
and $(1 \mapsto \phi(2), 2 \mapsto \phi(1))$,
and note that these form an orbit
(and similarly for the other edges).

\subsection{Decompositions}

Now we set up the general framework for decompositions.

\begin{defn} \label{def:vsys} 
Let $\mc{A}$ be a set of $R$-complexes;
we call $\mc{A}$ an $R$-complex family.
If each $A \in \mc{A}$ is a copy of $\Ss^\le$
we call $\mc{A}$ a $\Ss^\le$-family.
For $r \in \mb{N}$ we write 
$A_r = \bigcup \{ A_B: B \in \tbinom{R}{r} \}$
and $\mc{A}_r = \cup_{A \in \mc{A}} A_r$.

Let $\Phi$ be an $R$-complex.
We let $A(\Phi)$ denote the set of $\Phi$-embeddings of $A$.
We let $A(\Phi)^\le$ denote the $V(A)$-complex
where each $A(\Phi)^\le_F$ for $F \sub V(A)$
is the set of $\Phi$-embeddings of $A[F]$.

We let $\mc{A}(\Phi)^\le$ denote the 
$V(A)$-complex family $(A(\Phi)^\le: A \in \mc{A})$.

Let $\gG \in \GG^{\mc{A}_r}$ for some abelian group $\GG$.

For $\phi \in A(\Phi)^\le$ with $A \in \mc{A}$
we define $\gG(\phi) \in \GG^{\Phi_r}$ 
by $\gG(\phi)_{\phi \circ \tT} = \gG_\tT$ for $\tT \in A_r$
(zero otherwise).
For $\phi \in A(\Phi)$ we call $\gG(\phi)$ a $\gG$-molecule.
We let $\gG(\Phi)$ be the set of $\gG$-molecules.

Given $\Psi \in \mb{Z}^{\mc{A}(\Phi)}$ 
we define $\pl \Psi = \pl^\gG \Psi
= \sum_\phi \Psi_\phi \gG(\phi) \in \GG^{\Phi_r}$.
We also call $\Psi$ an integral 
$\gG(\Phi)$-decomposition of $\pl \Psi$
and call $\bgen{\gG(\Phi)}$ the decomposition lattice.
If furthermore $\Psi \in \{0,1\}^{\mc{A}(\Phi)}$
(i.e.\ $\Psi \sub \mc{A}(\Phi)$) we call $\Psi$ 
a $\gG(\Phi)$-decomposition.
\end{defn}

In our example, $\mc{A} = \{A\}$ consists of a single
copy of the $[3]$-complex $\Ss^\le$ on $[3]$,
which is identical with $H$ as above, i.e.\ the complete 
$(\{1,2\},3)$-partite $[3]$-complex wrt $(\{1,2\},3)$.
We let $\GG = \mb{Z}^3$ and denote 
the standard basis by $e_{12}$, $e_{13}$, $e_{23}$, 
which we think of as the colours red, blue and green.
We define $\gG \in \GG^{\mc{A}_2}$ 
by $\gG_\tT = e_{Im(\tT)}$.
The constituent parts of our decompositions
are $\gG$-molecules $\gG(\phi)$, which
encode rainbow triangles in $\Phi$:
we have $\phi \in A(\Phi)$
(which can be identified\footnote{
This identification is convenient but perhaps 
potentially confusing: depending on the context,
we may identify the domain of $\phi$ with
either the domain or the range of maps in $\Ss$.}
with $\Phi_3$), i.e.\ $\phi \circ \tT \in \Phi$ 
for all $\tT \in A = \Ss^\le$,
and e.g.\ the blue edge $\phi(1)\phi(3)$
is encoded by the coordinates
$\gG(\phi)_{\phi \circ \tT} = \gG_\tT = e_{13}$
for $\tT \in A_2$ with $Im(\tT) = \{1,3\}$,
i.e.\ $\tT = (1 \mapsto 1, 3 \mapsto 3)$
and $\tT = (2 \mapsto 1, 3 \mapsto 3)$.
We encode any coloured graph $G$ by
$G^* \in (\mb{Z}^3)^{\Phi_2}$ defined by
$G^*_\psi = e_{12}$ if $Im(\psi)$ is a red edge,
$G^*_\psi = e_{13}$ if $Im(\psi)$ is a blue edge,
$G^*_\psi = e_{23}$ if $Im(\psi)$ is a green edge.
Then a $\gG(\Phi)$-decomposition of $G^*$
encodes a rainbow triangle decomposition of $G$.

Now we formalise in general the objects (atoms) 
that are being decomposed into molecules.

\begin{defn} \label{def:atom} (atoms)
For any $\phi \in \mc{A}(\Phi)$ and $O \in \Phi_r/\Ss$
such that $\gG(\phi)^O \ne 0$ we call
$\gG(\phi)^O$ a $\gG$-atom at $O$.
We write $\gG[O]$ for the set of $\gG$-atoms at $O$.
We say $\gG$ is elementary if
all $\gG$-atoms are linearly independent.
We define a partial order $\le_\gG$
on $\GG^{\Phi_r}$ where $H \le_\gG G$
iff $G-H$ can be expressed as the sum
of a multiset of $\gG$-atoms.
\end{defn}

In our example, atoms represent coloured edges.
To see this, consider again the encoding of the
blue edge $\phi(1)\phi(3)$ described above.
The relevant orbit $O \in \Phi_2/\Ss$
consists of the two labelled edges
$(1 \mapsto \phi(1), 3 \mapsto \phi(3))$
and $(2 \mapsto \phi(1), 3 \mapsto \phi(3))$,
and the relevant $\gG$-atom at $O$ 
is $\gG(\phi)^O$ which is a vector supported on $O$
with both coordinates equal to $e_{13}$.
There are two other $\gG$-atoms at $O$,
which are vectors supported on $O$
with both coordinates equal to $e_{12}$ (meaning red edge),
or both coordinates equal to $e_{23}$ (meaning green edge).
Thus $\gG$ is elementary, which is an important assumption
in our main theorem, ensuring that our decomposition
problems do not exhibit arithmetic peculiarities
(as seen e.g.\ in the Frobenius coin problem).

\subsection{Lattices}

We conclude with a characterisation of
the decomposition lattice $\bgen{\gG(\Phi)}$,
with conditions that are somewhat analogous
to the degree-based divisibility conditions
considered above, but also account for the
labels on the edges and the orbits of the group action.
Throughout we let $\Phi$ be a $\Ss$-adapted $[q]$-complex
for some $\Ss \le S_q$, let $\mc{A}$ be a $\Ss^\le$-family
and $\gG \in \GG^{\mc{A}_r}$. 

\begin{defn} \label{def:L}
For $J \in \GG^{\Phi_r}$ we define
$J^\sharp \in (\GG^Q)^\Phi$ by\footnote{
The notation $\psi' \sub \psi$ means
that $\psi'$ is a restriction of $\psi$.} 
$(J^\sharp_{\psi'})_B = 
\sum \{ J_\psi: \psi' \sub \psi \in \Phi_B \}$
for $B \in Q := \tbinom{[q]}{r}$, $\psi' \in \Phi$.
We define $\gG^\sharp \in (\GG^Q)^{\cup \mc{A}}$
by $(\gG^\sharp_{\tT'})_B = 
\sum \{ \gG_\tT: \tT' \sub \tT \in A_B \}$
for $B \in Q$, $\tT' \in A \in \mc{A}$.
We let $\mc{L}_\gG(\Phi)$ be the set of all
$J \in \GG^{\Phi_r}$ such that
$(J^\sharp)^O \in \sgen{\gG^\sharp[O]}$
for any $O \in \Phi/\Ss$. 
\end{defn}

We illustrate Definition \ref{def:L} with our running example.
We start with the orbit $O = \{\es\}$,
where $\es$ denotes the unique function
with domain $\es$ (also denoting the empty set).
Recall that we encode our coloured graph $G$
by $G^* \in (\mb{Z}^3)^{\Phi_2}$ and 
write $G^{ij}$ for the edges of $G$ 
with colour corresponding to $ij$.
Then $((G^*)^\sharp_\es)_{ij} =
\sum_{\psi \in \Phi_{ij}} G^*_\psi$
equals $2|G^{12}|e_{12}$ if $ij=12$
or $|G^{13}|e_{13}+|G^{23}|e_{23}$ otherwise.
Similarly, $(\gG^\sharp_\es)_{ij} 
= \sum_{\tT \in \Ss^\le_{ij}} \gG_\tT$
equals $2e_{12}$ if $ij=12$ 
or $e_{13}+e_{23}$ otherwise.
The $0$-divisibility condition is that
$(2|G^{12}|e_{12}, |G^{13}|e_{13}+|G^{23}|e_{23},
|G^{13}|e_{13}+|G^{23}|e_{23} )$ 
is an integer multiple of
$(2e_{12},e_{13}+e_{23},e_{13}+e_{23})$, i.e.\
$G$ has an equal number of edges of each colour.

Next consider the $1$-divisibility condition
for any orbit $O = \{(1 \to x), (2 \to x)\}$ with $x \in V_1$. 
For $i,i' \in [2]$, $j \ne i$ we have
$((G^*)^\sharp_{i \to x})_{ij} 
= \sum \{ G^*_\psi: \psi \in \Phi_{ij}, \psi(i) = x\} $,
which equals $|G^{12}(x)|e_{12}$ if $j \in [2]$
or $|G^{13}(x)|e_{13}+|G^{23}(x)|e_{23}$ if $j=3$. 
Also, $(\gG^\sharp(i' \to x)_{i \to x})_{ij}
= (\gG^\sharp_{i \to i'})_{ij}
= \sum \{ \gG_\tT: \tT \in \Ss^\le_{ij}, \tT(i)=i' \}$,
which equals $e_{12}$ if $j \in [2]$
or $e_{i'3}$ if $j=3$. Thus we need
$(|G^{12}(x)|e_{12}, |G^{13}(x)|e_{13}+|G^{23}(x)|e_{23},
|G^{13}(x)|e_{13}+|G^{23}(x)|e_{23})$ to lie in the group
generated by $(e_{12},e_{13},0)$, $(e_{12},e_{23},0)$,
$(e_{12},0,e_{13})$ and $(e_{12},0,e_{23})$,
which holds iff $|G(x) \cap V_1|=|G(x) \cap V_2|$,
i.e.\ each $x \in V_1$ has equal degrees in $V_1$ and in $V_2$.

The other $1$-divisibility conditions are 
for orbits $O = \{ 3 \to x\}$ with $x \in V_2$.
For $i \in [2]$ we have
$((G^*)^\sharp_{3 \to x})_{i3} 
= \sum \{ G^*_\psi: \psi \in \Phi_{i3}, \psi(3) = x\} 
= |G^{13}(x)|e_{13}+|G^{23}(x)|e_{23}$ and
$(\gG^\sharp(3 \to x)_{3 \to x})_{i3}
= (\gG^\sharp_{3 \to 3})_{i3}
= \sum \{ \gG_\tT: \tT \in \Ss^\le_{i3}, \tT(3)=3 \}
= e_{13} + e_{23}$, so we need $|G^{13}(x)| = |G^{23}(x)|$,
i.e.\ each $x \in V_2$ has blue degree equal to green degree.
There are no further conditions, as the 
$2$-divisibility conditions hold trivially
(we leave this verification to the reader).

Returning to the general setting,
it is not hard to see $\bgen{\gG(\Phi)} \sub \mc{L}_\gG(\Phi)$.
The following result (\cite[Lemma 5.19]{K2}) shows that the 
converse inclusion holds under an extendability assumption on $\Phi$.

\begin{lemma} \label{lattice}
Let $\Ss \le S_q$, $\mc{A}$ be a $\Ss^\le$-family
and $\gG \in (\mb{Z}^D)^{\mc{A}_r}$. 
Let $\Phi$ be a $\Ss$-adapted 
$(\oO,s)$-extendable $[q]$-complex with $s=3r^2$, 
$n = |V(\Phi)| > n_0(q,D)$ large and $\oO > n^{-1/2}$.
Then $\bgen{\gG(\Phi)} = \mc{L}_\gG(\Phi)$.
\end{lemma}

\subsection{Types and regularity}

Next we will formulate our regularity assumption, 
which can be thought of as robust fractional decomposition.
First we give another notation for atoms. 

\begin{defn} \label{def:atom=}
For $\psi \in \Phi_B$ and $\tT \in \mc{A}_B$
we define $\gG[\psi]^\tT \in \GG^{\psi\Ss}$
by $\gG[\psi]^\tT_{\psi\sS} = \gG_{\tT\sS}$.
\end{defn}

We will illustrate the various notations in our example for the atom 
$\gG(\phi)^O$ representing a blue edge $\phi(1)\phi(3)$ as above.
In the notation of Definition \ref{def:vsys}, 
we write $\gG(\phi)^O=\gG(\phi')$
where $\phi'=\phi\mid_{\{1,3\}}$ has domain $\{1,3\}$,
so if $\tT \in A_2$ with $Im(\tT) \sub Dom(\phi')$
then $\tT = (1 \mapsto 1, 3 \mapsto 3)$
or $\tT = (2 \mapsto 1, 3 \mapsto 3)$.
In the notation of Definition \ref{def:atom=}, 
we write $\gG(\phi)^O=\gG[\phi']^\tT$
with $\tT = (1 \mapsto 1, 3 \mapsto 3)$,
as $\gG[\phi']^\tT$ is supported on
$\phi' = (1 \mapsto \phi(1), 3 \mapsto \phi(3))$
with value $\gG_\tT =  e_{13}$ and on
$\phi' \circ (12) = (2 \mapsto \phi(1), 3 \mapsto \phi(3))$
with value $\gG_{\tT \circ (12)} =  e_{13}$.
We also think of this notation as 
`an atom of type $\tT$ on $\psi$', 
where we define types in general as follows.

\begin{defn} \label{def:type} (types)
For $\tT \in \mc{A}_B$ we define 
$\gG^\tT \in \GG^{\Ss^B}$ by $\gG^\tT_\sS = \gG_{\tT\sS}$.

A type $t=[\tT]$ in $\gG$ is an equivalence class 
of the relation $\sim$ on any $\mc{A}_B$ 
with $B \in Q = \tbinom{[q]}{r}$
where $\tT \sim \tT'$ iff $\gG^\tT = \gG^{\tT'}$.
We write $T_B$ for the set of types in $\mc{A}_B$.

For $\tT \in t \in T_B$ and $\psi \in \Phi_B$ we write 
$\gG^t = \gG^\tT$ and $\gG[\psi]^t = \gG[\psi]^\tT$.

If $\gG^t=0$ call $t$ a zero type and write $t=0$.

If $\phi \in \mc{A}(\Phi)$ with
$\gG(\phi)^{\psi\Ss} = \gG[\psi]^t$
we write $t_\phi(\psi) = t$.
\end{defn}

To illustrate the preceding 
definition on the above example of $\gG[\phi']^\tT$
with $\tT = (1 \mapsto 1, 3 \mapsto 3)$
we think of $\{\tT\} \in T_{13}$ as the `blue edge' type with
$(\gG[\phi']^\tT_{\phi'}, \gG[\phi']^\tT_{\phi' \circ (12) })
= (\gG^\tT_{id}, \gG^\tT_{(12)}) 
= (\gG_{1 \mapsto 1, 3 \mapsto 3}, \gG_{2 \mapsto 1, 3 \mapsto 3}) 
= (e_{13}, e_{13})$. The possibility of a zero type is not relevant
to our example, as it allows for non-edges when
decomposing into copies of a non-complete graph.
The `red edge' type in $T_{12}$
is $\{ (1 \mapsto 1, 2 \mapsto 2), 
(1 \mapsto 2, 2 \mapsto 1) \}$,
as $(\gG^{1 \mapsto 1, 2 \mapsto 2}_{id}, 
\gG^{1 \mapsto 1, 2 \mapsto 2}_{(12)})
= (\gG_{1 \mapsto 1, 2 \mapsto 2}, 
\gG_{1 \mapsto 2, 2 \mapsto 1}) = (e_{12},e_{12})$
and $(\gG^{1 \mapsto 2, 2 \mapsto 1}_{id}, 
\gG^{1 \mapsto 2, 2 \mapsto 1}_{(12)})
= (\gG_{1 \mapsto 2, 2 \mapsto 1}, 
\gG_{1 \mapsto 1, 2 \mapsto 2}) = (e_{12},e_{12})$.

Now we formulate our regularity assumption.
The following definition can be roughly understood 
as saying that the vector $J$ can be approximated
by a non-negative linear combination of molecules,
where all molecules that can be used 
(in that $J$ contains all their atoms)
are used with comparable weights
(up constant factors).

\begin{defn} (regularity)

Suppose $\gG$ is elementary and $J \in (\mb{Z}^D)^{\Phi_r}$
with $J^O \in \sgen{\gG[O]}$ for all $O \in \Phi_r/\Ss$.
For $\psi \in \Phi_B$ with $|B|=r$ we define 
integers $J^t_\psi$ for all nonzero $t \in T_B$ by
$J^{\psi\Ss} = \sum_{0 \ne t \in T_B} J^t_\psi \gG[\psi]^t$.
Any choice of orbit representatives $\psi^O \in \Phi_{B^O}$ 
for each orbit $O \in \Phi_r/\Ss$ defines an atom decomposition
$J = \sum_{O \in \Phi_r/\Ss} \sum_{0 \ne t \in T_{B^O}} J^t_{\psi^O} \gG[\psi^O]^t$.

Let $\mc{A}(\Phi,J) = \{ \phi \in \mc{A}(\Phi): \gG(\phi) \le_\gG J \}$.
We say $J$ is $(\gG,c,\oO)$-regular (in $\Phi$) if there is 
$y \in [\oO n^{r-q},\oO^{-1} n^{r-q}]^{\mc{A}(\Phi,J)}$ such that for all 
$B \in Q$, $\psi \in \Phi_B$, $0 \ne t \in T_B$ we have 
\[\pl^t y_\psi := \sum \{  y_\phi: t_\phi(\psi)=t \} = (1 \pm c)J^t_\psi.\]
\end{defn}

For example, suppose $J=G^* \in (\mb{Z}^3)^{\Phi_2}$ 
encodes $G$ as above. An atom decomposition expresses
$J$ as a sum where each summand encodes a coloured edge of $G$
by some atom $\gG[\psi^O]^t$ as discussed above.
We have $\phi \in \mc{A}(\Phi,J)$ iff the molecule
$\gG(\phi)$ encodes a rainbow triangle in $G$.
Then $G^*$ is $(\gG,c,\oO)$-regular if we can assign
each rainbow triangle in $G$ a weight between
$\oO n^{-1}$ and $\oO^{-1} n^{-1}$ so that the total weight
of triangles on any edge is $1 \pm c$.

We require one further definition, used in the
extendability hypothesis of Theorem \ref{main} below.

\begin{defn} \label{Gatoms}
For $L \in \GG^{\Phi_r}$ we let
$\gG[L]=(\gG[L]^A: A \in \mc{A})$
where each $\gG[L]^A$ is the set of
$\psi \in A(\Phi)^\le_r$ 
such that $\gG(\psi) \le_\gG L$.
\end{defn}

In our example, the extendability hypothesis
says that for any $\Phi$-extension 
$E = (J,F,\phi)$ of rank $h$ there are many 
$\phi^+ \in X_E(\Phi)$ such that all edges 
$Im(\phi^+\psi)$ with $\psi \in J_2 \sm J[F]$
are edges of $G$ with the correct colour
(red if $\psi \in J_{12}$,
blue if $\psi \in J_{13}$,
green if $\psi \in J_{23}$).
We illustrate this for
extensions of some fixed rainbow triangle
to an octahedron of rainbow triangles
(recall $J=[3](2)$, $F=[3] \times \{1\}=[3]$
and let $J' = J_2 \sm J[F]$). 
If $(\Phi,\gG[G^*]^A)$ is $(\oO,2)$-extendable
we have at least $\oO |V_1|^2 |V_2|$ choices
of $\phi^+ \in X_{E,J'}(\Phi,\gG[G^*]^A)$.
For each $\psi \in J_2 \sm J[F]$
we have $\phi^+ \psi \in \gG[G^*]^A$,
i.e.\ $\gG(\phi^+\psi) \le_\gG G^*$.
For example, if $\psi \in J_{13}$ with
$\psi(1)=(1,1)$ and $\psi(3)=(3,2)$
then $\gG(\phi^+\psi)$ is the blue atom
at $Im(\phi^+\psi)$, i.e.\ 
the vector supported on the orbit 
with the two labelled edges
$(1 \mapsto \phi^+((1,1)), 3 \mapsto \phi^+((3,2)))$
and $(2 \mapsto \phi^+((1,1)), 3 \mapsto \phi^+((3,2)))$,
where both nonzero coordinates are $e_{13}$.
For this $\psi$, the condition 
$\gG(\phi^+\psi) \le_\gG G^*$ says that $G$ has
a blue edge at $\phi^+((1,1)) \phi^+((3,2))$.
As $\psi$ varies over $J_2$ we see that $Im(\phi^+)$
spans an octahedron of rainbow triangles.

Finally we can state the main result (Theorem 3.1) of \cite{K2}
(recall $h=2^{50q^3}$ and $\dD = 2^{-10^3 q^5}$).

\begin{theo} \label{main}
For any $q \ge r$ and $D$ there are $\oO_0$ and $n_0$ such that 
the following holds for $n > n_0$, 
$n^{-\dD}<\oO<\oO_0$ and $c \le \oO^{h^{20}}$.
Let $\mc{A}$ be a $\Ss^\le$-family with $\Ss \le S_q$.
Suppose $\gG \in (\mb{Z}^D)^{\mc{A}_r}$ is elementary.
Let $\Phi$ be a $\Ss$-adapted $[q]$-complex on $[n]$. 
Let $G \in \bgen{\gG(\Phi)}$ be $(\gG,c,\oO)$-regular 
in $\Phi$ such that $(\Phi,\gG[G]^A)$ 
is $(\oO,h)$-extendable for each $A \in \mc{A}$.
Then $G$ has a $\gG(\Phi)$-decomposition.
\end{theo}

\section{Coloured hypergraphs} \label{sec:col}

When can an edge-coloured graph be
decomposed into rainbow triangles?
In this section we illustrate the application of Theorem \ref{main}
to this question, and a hypergraph generalisation thereof.
We start by formulating the general problem
of decomposing an edge-coloured $r$-multigraph $G$
by an edge-coloured $r$-graph $H$.
For simplicity we assume that $H$ is simple
(one could allow multiple copies of edges in $H$
provided they have distinct colours, but not multiple
edges of a given colour, as then the associated $\gG$ 
in Definition \ref{def:col*2} below is not elementary). 

\begin{defn} \label{def:colhyp}
Suppose $H$ is an $r$-graph on $[q]$,
edge-coloured as $H = \cup_{d \in [D]} H^d$.
We identify $H$ with a vector $H \in (\mb{N}^D)^Q$,
where each $(H_f)_d = 1_{f \in H^d}$ (indicator function)
and $Q=\tbinom{[q]}{r}$.

Let $\Phi$ be an
$S_q$-adapted $[q]$-complex on $[n]$.
For $\phi \in \Phi_q$ we define 
$\phi(H) \in (\mb{N}^D)^{\Phi^\circ_r}$
by $\phi(H)_{\phi(f)} = H_f$.
Let $\mc{H}$ be an family of 
$[D]$-edge-coloured $r$-graphs on $[q]$.
Let $\mc{H}(\Phi) = \{ \phi(H):
 \phi \in \Phi_q, H \in \mc{H} \}$.

Let $G \in \mb{N}^{\Phi^\circ_r}$ be an $r$-multigraph 
$[D]$-edge-coloured as $G = \cup_{d \in [D]} G^d$,
identified with $G \in (\mb{N}^D)^{\Phi^\circ_r}$.
We call $\mc{H'} \sub \mc{H}(\Phi)$ with 
$\sum \mc{H'} = G$ an $H$-decomposition of $G$ in $\Phi$.
We call $\Psi \in \mb{Z}^{\mc{H}(\Phi)}$ with
$\sum_{H'} \Psi_{H'} H' = G$
an integral $H$-decomposition of $G$ in $\Phi$.
\end{defn}

Note that copies of $H$ in an integral $H$-decomposition of $G$
can use edges $e \in \Phi^\circ_r$ with $G_e=0$
or with the wrong colour, but all such terms must cancel.
Before considering the general setting of the previous definition,
we warm up by specialising to graphs $(r=2)$ and the case 
that $\Phi$ is the complete $[q]$-complex on $[n]$.
We formulate a typicality condition for coloured graphs
and a result on rainbow triangle decompositions 
analogous to that given in \cite{Kcount} for 
triangle decompositions of typical graphs.

\begin{defn} \label{def:coltyp}
Let $G$ be a $[D]$-edge-coloured graph on $[n]$.
For $\aA \in [D]$, the $\aA$-density of $G$
is $d(G^\aA) = |G^\aA| \tbinom{n}{2}^{-1}$.
The density of $G$ is $d(G) = |G| \tbinom{n}{2}^{-1}$.
The density vector of $G$ is $d(G)^* \in [0,1]^D$
with $d(G)^*_\aA = d(G^\aA)$.
Given vectors $\bm{x} \in [n]^t$ of vertices
and $\aB \in [D]^t$ of colours we define the $\aB$-degree  
$d^\aB_G(\bm{x})$ of $\bm{x}$ in $G$ as the number of vertices $y$ 
such that $x_i y \in G^{\aA_i}$ for all $i \in [t]$.

We say $G$ is $(c,h)$-typical if
$d^\aB_G(\bm{x}) = (1 \pm tc) n \prod_{i=1}^t d(G^{\aA_i})$
for any such $\bm{x}$ and $\aB$ with $t \le h$.
\end{defn}

\begin{theo} \label{rainbowtri}
Suppose $G$ is a tridivisible $(c,h)$-typical 
$[D]$-edge-coloured graph on $[n]$,
where $D \ge 4$, $n>n_0(D)$ is large, 
$h=2^{10^3}$, $\dD = 2^{-10^6}$, $c < c_0 d(G)^{h^{90}}$
where $c_0=c_0(D)$ is small, and each
$n^{-\dD/2h^3} < d(G^\aA) < (1/3 - n^{-\dD/2h^3}) d(G)$.
Then $G$ has a rainbow triangle decomposition. 
\end{theo}

Note that the tridivisibility condition
($G$ has all degrees even and $3 \mid e(G)$)
in Theorem \ref{rainbowtri} is necessary,
as if we ignore the colours then we obtain
a triangle decomposition of $G$;
it is perhaps surprising that the colours
do not impose any additional condition.
We will deduce Theorem \ref{rainbowtri} from a more
general result on typical $r$-multigraphs,
as in the following definition.

\begin{defn} \label{def:coltyp2}
Let $G$ be a $[D]$-edge-coloured $r$-multigraph on $[n]$.
For $\aA \in [D]$, the $\aA$-density of $G$
is $d(G^\aA) = |G^\aA| \tbinom{n}{r}^{-1}$.
The density of $G$ is $d(G) = |G| \tbinom{n}{r}^{-1}$.
The density vector of $G$ is $d(G)^* \in \mb{R}^D$
with $d(G)^*_\aA = d(G^\aA)$.

For $e \sub [n]$, the degree of $e$ in $G$ is $|G(e)|$;
the degree vector is $G(e)^* \in \mb{N}^D$ 
with $G(e)^*_\aA = |G^\aA(e)|$.

Given vectors
$\bm{f} \in \tbinom{[n]}{r-1}^t$ of $(r-1)$-sets
and $\aB \in [D]^t$ of colours we define 
the $\aB$-degree of $\bm{f}$ in $G$ as 
$d^\aB_G(\bm{f}) = \sum_{v \in [n]} 
\prod_{i=1}^t G^{\aA_i}_{f_i \cup \{v\}}$.

We say $G$ is $(c,h)$-typical if
$d^\aB_G(f) = (1 \pm tc) n \prod_{i=1}^t d(G^{\aA_i})$
for any such $\bm{f}$ and $\aB$ with $t \le h$.

Given a family $\mc{H}$ of $[D]$-edge-coloured $r$-graphs on $[q]$,
we say $G$ is $(b,c)$-balanced wrt $\mc{H}$ if 
there is $p \in [b,b^{-1}]^{\mc{H}}$ with
$d(G)^* = (1 \pm c) \sum_H p_H d(H)^*$.

We say $G$ is $\mc{H}$-divisible if each
$G(e)^* \in \sgen{ H(f)^*: f \in \tbinom{[q]}{|e|}, H \in \mc{H} }$.
\end{defn}

In the next lemma we show that in the case of rainbow triangles,
the conditions in Definition \ref{def:coltyp2}
follow from the assumptions of Theorem \ref{rainbowtri}.

\begin{lemma} \label{tris}
Let $\mc{H}$ be the family of all $[D]$-edge-coloured rainbow triangles
and $G$ be a $[D]$-edge-coloured graph on $[n]$, with $D \ge 4$. Then
\begin{enumerate}
\item $G$ is $\mc{H}$-divisible iff $G$ is tridivisible, and
\item If each $bD^2 < d(G^\aA) < d(G)/3 - bD^3$ 
then $G$ is $(b,0)$-balanced wrt $\mc{H}$.
\end{enumerate}
\end{lemma}

\nib{Proof.}
For (i), we need to know the integer span $Z(r,s)$ of the rows
of a matrix $M(r,s)$ whose rows are indexed by $\tbinom{[s]}{r}$
and columns by $[s]$, with $M(r,s)_{e,i} = 1_{i \in e}$.
It follows from \cite[Theorem 2]{W7} 
(and is not hard to show directly) that 
$Z(r,s) = \{ \bm{x} \in \mb{Z}^s: r \mid \sum_i x_i\}$ for $s>r$.
To apply this to the divisibility conditions,
first consider $G(\es)^* = (|G^1|,\dots,|G^D|)$ and note that 
$H(\es)^* =  (|H^1|,\dots,|H^D|)$ for $H \in \mc{H}$
are the rows of $M(3,D)$. We have
$G(\es)^* \in \sgen{ H(\es)^*: H \in \mc{H}}$
iff $3 \mid \sum_\aA |G^\aA| = |G|$.
Next, for any $v \in [n]$ we have
$G(v)^* = (|G^1(v)|,\dots,|G^D(v)|)$.
As $H(x)^* =  (|H^1(x)|,\dots,|H^D(x)|)$ 
for $x \in [q]$, $H \in \mc{H}$ are the rows of $M(2,D)$
we have $G(v)^* \in \sgen{ H(x)^*: x \in [q], H \in \mc{H}}$
iff $2 \mid \sum_\aA |G^\aA(v)| = |G(v)|$.
Finally, for any $uv \in \tbinom{[n]}{2}$ we have
$G(uv)^* = (G^1_{uv},\dots,G^D_{uv})$ and
$H(xy)^*$ for $xy \in \tbinom{[q]}{2}$, $H \in \mc{H}$
is the standard basis, so the $2$-divisibility condition is trivial.
Thus $G$ is $\mc{H}$-divisible iff $G$ is tridivisible.

For (ii), we note that the set of density vectors $d(H)^*$ 
for $H \in \mc{H}$ consists of all probability distributions
on $[D]$ with $3$ coordinates equal to $1/3$ and the rest zero.
By \cite[Theorem 46]{HLP}, any probability distribution $\bm{x}$ 
on $[D]$ is a convex combination of the vectors $d(H)^*$ 
iff $x_\aA \le 1/3$ for all $\aA \in [D]$. Thus for
any $\bm{x} \in [0,1]^D$ with each $3x_{\aA'} \le \sum_\aA x_\aA \le 1$
there is some $\bm{p} \in [0,1]^{\mc{H}}$ with $\bm{x} = \sum_H p_H d(H)^*$ 
and $\sum_H p_H = \sum_\aA x_\aA$. We apply this to
$\bm{x} = d(G)^* - b \sum_H d(H)^*$, noting that
$\sum_\aA x_\aA = d(G) - b \tbinom{D}{3}$ and
each $0 \le x_\aA = d(G^\aA) - \tfrac{b}{3} \tbinom{D-1}{2} 
\le \tfrac{1}{3} \sum_\aA x_\aA$.
Then $\bm{p}' = \bm{p}+b\bm{1} \in [b,b^{-1}]^{\mc{H}}$ has
$d(G)^* = \sum_H p'_H d(H)^*$. \qed

\medskip

Next we consider how to encode decompositions
of coloured multigraphs in the labelled edge setting
of Theorem \ref{main}; this is similar to the running
example used in the previous section.

\begin{defn} \label{def:col*}
Given a set $e$ of size $r$, we write
$e^{r \to q}$ for the set of all $\pi^{-1}$
where $\pi:e \to [q]$ is injective.
Given a $[D]$-edge-coloured $r$-multigraph 
$G = (G^d: d \in [D])$ we define
$G^{r \to q} = ((G^{r \to q})^d: d \in [D])$
where each $(G^{r \to q})^d$ is the (disjoint) union of
all $e^{r \to q}$ with $e \in G^d$.
\end{defn}

\begin{lemma} \label{equivcol*}
Let $H$ and $G$ be $[D]$-edge-coloured $r$-multigraphs, 
$H^* = H^{r \to q}$ and $G^* = G^{r \to q}$.
Then an (integral) $H$-decomposition of $G$ is equivalent to
an (integral) $H^*$-decomposition of $G^*$. 
\end{lemma}

\nib{Proof.}
We associate any $H$-decomposition $\mc{D}$ of $G$
with an $H^*$-decomposition $\mc{D}^*$ of $G^*$, associating 
each $\phi(H) \in \mc{D}$ with 
$\phi H^* := \{ \phi \circ \tT: \tT \in H^* \} \in \mc{D}^*$.
Then $e \in \phi(H^d)$ iff $e^{r \to q} \sub \phi H^{*d}$,
as if $e = \phi(f)$ for some $f \in H^d$
and $\pi^{-1} \in e^{r \to q}$ then $\pi^{-1} = \phi \tT$, 
where $\tT = \phi^{-1} \pi^{-1} \in H^{*d}$, and conversely. 
The same proof applies to integral decompositions
(defined in Definition \ref{def:colhyp}). \qed

\begin{defn} \label{def:col*2}
Given a family $\mc{H}$ of $[D]$-edge-coloured $r$-graphs on $[q]$,
let $\mc{A} =\mc{A}^{\mc{H}}=\{A^H: H \in \mc{H}\}$ with each $A^H=S_q^\le$
and $\gG = \gG^{\mc{H}} \in (\mb{Z}^D)^{\mc{A}_r}$ with $\gG_\tT = e_d$
(standard basis vector)
if $\tT \in A^H_r$, $H \in \mc{H}$, $d \in [D]$
with $Im(\tT) \in H^d$ or $\gG_\tT=0$ otherwise.
\end{defn}

\begin{lemma} \label{coldiv}
With notation as in Definitions 
\ref{def:colhyp}, \ref{def:col*} and \ref{def:col*2},
an (integral) $\mc{H}$-decomposition of $G$ in $\Phi$ is equivalent 
to an (integral) $\gG(\Phi)$-decomposition of $G^*$.

Furthermore, if $\Phi$ is $(\oO,s)$-extendable 
with $s=3r^2$, $\oO>n^{-1/2}$ and $n>n_0(q)$ large
then $G$ has an integral $\mc{H}$-decomposition in $\Phi_q$
iff $G$ is $\mc{H}$-divisible.
\end{lemma}

\nib{Proof.}
For the first statement, the same argument
as in Lemma \ref{equivcol*} shows that 
an $\mc{H}$-decomposition of $G$ in $\Phi$ is equivalent to
an $\mc{H}^*$-decomposition of $G^*$ in $\Phi$
(where $\mc{H}^* = \{H^*: H \in \mc{H}\}$),
i.e.\ some $\mc{D} \sub \mc{H}^*(\Phi)
= \{ \phi H^*: H \in \mc{H}, \phi \in \Phi_q \}$
with $\sum \mc{D} = G^* \in (\mb{N}^D)^{\Phi_r}$.
We can also view $\mc{D}$ as a $\gG(\Phi)$-decomposition 
of $G^*$ by identifying each $\phi H^* \in \mc{D}$ 
with the molecule $\gG(\phi)$ where $\phi \in A^H(\Phi)$:
indeed, if $\phi \pi^{-1} \in \phi H^{*d}$ with $d \in [D]$, 
where $e \in H^d$ and $\pi:e \to [q]$ is injective, then
$\gG(\phi)_{\phi \pi^{-1}} = \gG_{\pi^{-1}} = e_d$.
This proves the equivalence for decompositions,
and the same argument applies to integral decompositions.

For the second statement, by Lemma \ref{lattice} we have
$\bgen{\gG(\Phi)} = \mc{L}_\gG(\Phi)$.
By Definition \ref{def:L} we need to show
that $G$ is $\mc{H}$-divisible iff
$((G^*)^\sharp)^O \in \sgen{\gG^\sharp[O]}$
for any $O \in \Phi/S_q$. 

Fix any $O \in \Phi/S_q$,
write $e=Im(O) \in \Phi^\circ$ and $i=|e|$.
Then $((G^*)^\sharp)^O \in ((\mb{Z}^D)^Q)^O 
= (\mb{Z}^D)^{Q \times O}$
is a vector supported on the coordinates
$(B,\psi')$ with $B' \sub B \in Q$ and $\psi' \in O \cap \Phi_{B'}$
with each $((G^*)^\sharp_{\psi'})_B) = 
\sum \{ G^*_\psi: \psi' \sub \psi \in \Phi_B \} 
= (r-i)! G(e)^* \in \mb{N}^D$.

Also, $\sgen{\gG^\sharp[O]}$ is generated by $\gG^\sharp$-atoms
$\gG^\sharp(\ups)$ at $O$, each of which is
supported on the same coordinates $(B,\psi')$
as $((G^*)^\sharp)^O$, with each
$(\gG^\sharp(\ups)_{\psi'})_B)$ equal 
to some $(r-i)! H(f)^*$ with $f \in \tbinom{[q]}{|e|}$, $H \in \mc{H}$.
The lemma follows. \qed

\medskip

Now we state our theorem on decompositions
of typical coloured $r$-multigraphs.
By Lemma \ref{tris} it implies Theorem \ref{rainbowtri}.
We will deduce it from Theorem \ref{colHdecomp:ext} below.

\begin{theo} \label{coltypmulti}
Let $\mc{H}$ be a family of 
$[D]$-edge-coloured $r$-graphs on $[q]$.
Suppose $G$ is a $(c,h^q)$-typical 
$[D]$-edge-coloured $r$-multigraph on $[n]$
with all $G^d_e < b^{-1}$
that is $(b,c)$-balanced wrt $\mc{H}$,
where $n>n_0(q,D)$ is large,
$d(G) > b := n^{-\dD/h^q}$, 
$c < c_0 d(G)^{h^{30q}}$
and $c_0=c_0(q)$ is small.
Then $G$ has an $\mc{H}$-decomposition
iff $G$ is $\mc{H}$-divisible.
\end{theo}

The next definition formulates the extendability and regularity
conditions for coloured hypergraph decompositions;
we will see below that they both follow from typicality.
We remark that the extendability condition is stronger than
simply requiring that each $(\Phi,G^i)$ is extendable
(it is roughly equivalent to certain lower bounds on 
degree vectors $d^\gG_G(x)$ as in Definition \ref{def:coltyp}).

\begin{defn} \label{def:colhyp2}
With notation as in Definition \ref{def:colhyp},
we say $G \in (\mb{N}^D)^{\Phi^\circ_r}$ 
is $(\mc{H},c,\oO)$-regular in $\Phi$ if there are
$y^H_\phi \in [\oO n^{r-q},\oO^{-1} n^{r-q}]$ for each 
$H \in \mc{H}$, $\phi \in \Phi_q$ with $\phi(H) \le G$
(coordinate-wise) so that 
$\sum \{ y^H_\phi \phi(H) \} = (1 \pm c)G$ 
(sum over all valid $(H,\phi)$, 
approximation coordinate-wise).

We say that $(\Phi,G)$ is $(\oO,h)$-extendable
if $(\Phi,G')$ is $(\oO,h)$-extendable,
where $G' = (G^1,\dots,G^D)$.
\end{defn}

The next theorem shows extendability and regularity suffice
for the equivalence of decomposition and integral decomposition.
For wider applicability we formulate it in the setting of
exactly adapted complexes, as in the following definition,
which allows for an $S_q$-adapted $[q]$-complex
(such as the complete $[q]$-complex, suppressed in the
statement of Theorem \ref{coltypmulti}),
or a generalised partite complex, which is
exactly $\Ss$-adapted for some subgroup $\Ss$ of $S_q$
(such as that in the running example of the previous section).

\begin{defn}
We say that an $R$-complex $\Phi$ is exactly $\Ss$-adapted 
if whenever $\phi \in \Phi_B$ and $\tau \in Bij(B',B)$ 
(set of bijections from $B'$ to $B$) we have
$\phi \circ \tau\in \Phi_{B'}$ iff $\sS \in \Ss^B_{B'}$.

We say $\Phi$ is exactly adapted if 
$\Phi$ is exactly $\Ss$-adapted for some $\Ss$.
\end{defn}

\begin{theo} \label{colHdecomp:ext}
Let $\mc{H}$ be an family of 
$[D]$-edge-coloured $r$-graphs on $[q]$.
Let $\Phi$ be an $(\oO,h)$-extendable 
exactly adapted $[q]$-complex on $[n]$
where $n>n_0(q,D)$ is large,
$n^{-\dD}<\oO<\oO_0(q,D)$ is small and $c=\oO^{h^{20}}$.
Suppose $G \in (\mb{N}^D)^{\Phi^\circ_r}$ 
is $(\mc{H},c,\oO)$-regular in $\Phi$
and $(\Phi,G)$ is $(\oO,h)$-extendable.
Then $G$ has an $\mc{H}$-decomposition in $\Phi_q$
iff $G$ has an integral $\mc{H}$-decomposition in $\Phi_q$.
\end{theo}

\nib{Proof.} 
By Lemma \ref{coldiv}, it is equivalent to consider 
$\gG(\Phi)$-decompositions of $G^*$, with notation as in 
Definitions \ref{def:col*} and \ref{def:col*2}.
There are $D+1$ types in $\gG$ for each $B \in Q$:
the colour $d$ type 
$\{\tT \in A^H_B: Im(\tT) \in H^d, H \in \mc{H} \}$
for each $d \in [D]$, and the nonedge type 
$\{\tT \in A^H_B: Im(\tT) \notin H \in \mc{H} \}$.
Each $\gG^\tT$ is $e_d$ in all coordinates
for $\tT$ in a colour $d$ type
or $0$ in all coordinates
for $\tT$ in a nonedge type,
so $\gG$ is elementary.
The atom decomposition of $G^*$ is
$G^* = \sum_{f \in \Phi^\circ_r} \sum_{d \in [D]} 
(G_f)_d f^d$, where $f^d_\psi = e_d f^{r \to q}$. 

As $G$ is $(\mc{H},c,\oO)$-regular in $\Phi$ we have
$\sum \{ y^H_\phi \phi(H) \} = (1 \pm c)G$ for some
$y^H_\phi \in [\oO n^{r-q},\oO^{-1} n^{r-q}]$ for each 
$H \in \mc{H}$, $\phi \in \Phi_q$ with $\phi(H) \le G$.
As in the proof of the first part of Lemma \ref{coldiv},
we can identify any such $\phi(H)$ with $\phi H^* \le G^*$,
and so (regarding $\phi \in A^H(\Phi)$) 
with $\gG(\phi) \le_\gG G^*$, so $\phi \in \mc{A}(\Phi,G^*)$. 
Let  $y_\phi = y^H_\phi$ for $\phi \in A^H(\Phi)$.
For any $B \in Q$, $\psi \in \Phi_B$, $d \in [D]$,
writing $t_d \in T_B$ for the colour $d$ type,
$\pl^{t_d} y_\psi = \sum \{ y_\phi : t_\phi(\psi)=t_d \} 
= \sum \{ y^H_\phi: Im(\psi) \in \phi(H^d), H \in \mc{H} \}
= (1 \pm c)(G^*)^{t_d}_\psi$,
so $G^*$ is $(\gG,c,\oO)$-regular.

To apply Theorem \ref{main}, it remains to show that 
each $(\Phi,\gG[G^*]^H)$ is $(\oO,h)$-extendable.
If $B \notin H$ then $\gG[G^*]^H_B = \Phi_B$
and if $B \in H^d$ for $d \in [D]$ then 
$\gG[G^*]^H_B = \{ \psi \in \Phi_B : Im(\psi) \in G^d \}$.
Consider any $\Phi$-extension $E=(J,F,\phi)$ of rank $h$
and $J' \sub J_r \sm J[F]$.
Let $J'' = (J^d: d \in [D])$ with each
$J^d = \bigcup \{ J'_B: B \in H^d \}$.
As $(\Phi,G)$ is $(\oO,h)$-extendable
we have $X_{E,J''}(\Phi,G) > \oO n^{v_E}$.
Consider any $\phi^+ \in X_{E,J''}(\Phi,G)$.
For any $\psi \in J^d$ we have $\phi^+\psi \in \Phi$ 
and $Im(\phi^+\psi) \in G^d$, so $\phi^+\psi \in \gG[G^*]^H$.
Thus $\phi^+ \in X_{E,J'}(\Phi,\gG[G^*]^H)$,
so $(\Phi,\gG[G^*]^H)$ is $(\oO,h)$-extendable. \qed

\medskip

Now we show that the extendability and regularity
conditions follow from typicality, thus deducing
our decomposition result for typical coloured $r$-multigraphs.

\medskip

\nib{Proof of Theorem \ref{coltypmulti}.}
Suppose $G$ is an $\mc{H}$-divisible $(c,h^q)$-typical 
$[D]$-edge-coloured $r$-multigraph on $[n]$
that is $(b,c)$-balanced wrt $\mc{H}$,
where $n>n_0(q,D)$ is large,
$d(G) > b := 2n^{-\dD/h^q}$, 
$c < c_0 d(G)^{h^{30q}}$
and $c_0=c_0(q)$ is small.
We need to show that $G$ has an $\mc{H}$-decomposition.

Let $\Phi$ be the complete $[q]$-complex on $[n]$.
By Lemma \ref{coldiv} and $\mc{H}$-divisibility, 
$G$ has an integral $\mc{H}$-decomposition in $\Phi_q$.
Let $\bm{p} \in [b,b^{-1}]^{\mc{H}}$ with
$d(G)^* = (1 \pm c) \sum_H p_H d(H)^*$. We can assume
each colour $\aA \in [D]$ is used at least once by $\mc{H}$,
so $d(G^\aA) \ge b/2Q$, where $Q=\tbinom{q}{r}$.
To apply Theorem \ref{colHdecomp:ext},
it remains to check extendability and regularity.

We claim that $(\Phi,G)$ is $(\oO,h)$-extendable with $\oO > n^{-\dD}$.
To see this, consider any $\Phi$-extension 
$E = (J,F,\psi)$ with $J \sub [q](h)$
and $J' = (J^d: d \in [D])$ for some
mutually disjoint $J^d \sub J_r \sm J[F]$.
Let $V(J) \sm F = \{x_1,\dots,x_{v_E}\}$.
For $i \in [v_E]$ we list the neighbourhood 
of $x_i$ as $f^i = (f^i_1,\dots,f^i_{t_i})$ 
and let $\aB^i \in [D]^{[t_i]}$ be such that
each $f^i_j \cup \{x_i\}$ has colour $\aA^i_j$.
Then the number of choices for $x_i$
(weighted by edge-multiplicities)
given any previous choices $\phi'\mid_{\{x_j:j<i\}}$ 
is $d^{\aB^i}_G(\phi'(f^i)) 
= (1 \pm t_i c) n \prod_{j=1}^{t_i} d(G^{\aA^i_j})$.
As each $d(G^d) > b/2Q$ with $b = n^{-\dD/h^q}$,
we deduce 
\[X_{E,J'}(\Phi,G) = \sum_{\phi \in X_E(\Phi)} 
\prod_{d \in [D]} \prod_{f \in J^d} G^d_{\phi(f)}
> n^{v_E - \dD}.\]

For regularity, taking $E = (J,F,\psi)$ as above
with $J=[q](1)$, $J' = (H^d: d \in [D])$,
$F=f \in H^\aA$ with $H \in \mc{H}$, $\aA \in [D]$,
and $\psi \in Bij(f,e)$ with $e \in G^\aA$, we obtain
\[X_{E,J'}(\Phi,G) = (1 \pm Qc) d(G^\aA)^{-1}
n^{q-r} \prod_{d \in [D]} d(G^d)^{|H^d|}.\]
Let $Z = n^{q-r} \prod_{d \in [D]} d(G^d)^{|H^d|}$ and
$y_\phi = p_H (q)_r^{-1} Z^{-1} \prod_{d \in [D]} 
\prod_{f \in H^d} G^d_{\phi(f)}$
for each $\phi \in A^H(\Phi)$, $H \in \mc{H}$.
Then each such $y_\phi \in [\oO n^{r-q}, \oO^{-1} n^{r-q}]$,
as all $d(G^\dD) > b/2Q$, $p_H < b^{-1}$ and $G^d_{\phi(f)} < b^{-1}$.
Letting $f$ vary over $H^\aA$, we have
\begin{align*}
& \sum_H \sum_\phi y_\phi (\phi(H)_e)_\aA
= \sum_H p_H r! (q)_r^{-1} \sum_{f \in H^\aA} Z^{-1} 
\sum_{\phi \in X_E(\Phi)} \prod_{d \in [D]} 
\prod_{f \in H^d} G^d_{\phi(f)} \\
& = \sum_H p_H Q^{-1} \sum_{f \in H^\aA}
(1 \pm 2Q c) d(G^\aA)^{-1} G^\aA_e 
= (1 \pm q^r c) G^\aA_e.
\end{align*}
Thus $G$ is $(\mc{H},q^r c,\oO)$-regular in $\Phi$. \qed

\medskip

We conclude with a theorem on coloured 
generalised partite decompositions,
which can be used (we omit the details)
to obtain a common generalisation
of Theorems \ref{HPdecomp:typ} and \ref{coltypmulti}.

\begin{defn} \label{def:colblowup}
Let $\mc{H}$ be a family of 
$[D]$-edge-coloured $r$-graphs on $[q]$ and 
$\mc{P}=(P_1,\dots,P_t)$ be a partition of $[q]$.
Let $I^d = \{\iB \in \mb{N}^t: \cup_H H^d_\iB \ne \es\}$ 
and $I = \cup_d I^d$.

Let $\Ss$ be the group of all $\sS \in S_q$ 
with all $\sS(P_i)=P_i$.
Let $\Phi$ be an exactly $\Ss$-adapted $[q]$-complex
with parts $\mc{P}'=(P'_1,\dots,P'_t)$, where each 
$P'_i = \{ \psi(j): j \in P_i,\ \psi \in \Phi_{\{j\}} \}$.

Let $G \in (\mb{N}^D)^{\Phi^\circ_r}$. 
We call $G$ an $(\mc{H},\mc{P})$-blowup 
if $G^d_\iB \ne \es \Ra \iB \in I^d$. 

For $e \sub [n]$, $f \sub [q]$ 
we define $G(e)^*, H(f)^* \in (\mb{N}^D)^I$
by $(G(e)^*_\iB)_d=|G^d_\iB(e)|$, $(H(f)^*_\iB)_d=|H^d_\iB(f)|$.
We say $G$ is $(\mc{H},\mc{P})$-divisible if each
$G(e)^* \in \sgen{ H(f)^*: f \in \tbinom{[q]}{|e|}, H \in \mc{H} }$.
\end{defn}

In the following extendability hypothesis
we consider $G^d_\iB$ undefined for $\iB \notin I(H^d)$.

\begin{theo} \label{colpart}
With notation as in Definition \ref{def:colblowup},
suppose $n/h \le |P'_i| \le n$ with $n > n_0(q,D)$,
$G$ is an $(\mc{H},\mc{P})$-divisible $(\mc{H},\mc{P})$-blowup,
$G$ is $(\mc{H},c,\oO)$-regular in $\Phi$,
and $(\Phi,G)$ is $(\oO,h)$-extendable,
where $n^{-\dD}<\oO<\oO_0(q,D)$ and $c=\oO^{h^{20}}$.
Then $G$ has a $\mc{P}$-partite $\mc{H}$-decomposition.
\end{theo}

\nib{Proof.}
By Theorem \ref{colHdecomp:ext} it suffices to show that 
$G$ has an integral $\mc{H}$-decomposition in $\Phi_q$,
i.e.\ $G^* \in \bgen{\gG(\Phi)} = \mc{L}_\gG(\Phi)$ 
(by Lemmas \ref{coldiv} and \ref{lattice}).
Consider any $\iB \in I$
and $\iB' \in \mb{N}^t$ with all $i'_j \le i_j$.
Let $m^\iB_{\iB'} = \prod_{j \in [t]} (i_j-i'_j)!$.
For any $B' \sub B \in Q$ 
with $i_{\mc{P}}(B')=\iB'$ and  $i_{\mc{P}}(B)=\iB$
and $\psi' \in \Phi_{B'}$ with $Im(\psi')=e$ 
we have $((G^*)^\sharp_{\psi'})_B) = 
\sum \{ G^*_\psi: \psi' \sub \psi \in \Phi_B \} 
= m^\iB_{\iB'} G_\iB(e)^* \in \mb{N}^D$.
Writing $O=\psi'\Ss$, for any $\psi \in O$ we have 
$((G^*)^\sharp_\psi)_B) = m^\iB_{\iB'} G_\iB(e)^*$.
Thus we obtain $((G^*)^\sharp)^O$ from $G(e)^*$
by copying coordinates and multiplying all copies
of each $\iB$-coordinate by $m^\iB_{\iB'}$.
Similarly, for any 
$H \in \mc{H}$, $\tT' \in A^H_{B'}$, $f=Im(\tT')$ 
we have $(\gG^\sharp_{\tT'})_B) = 
\sum \{ \gG_\tT: \tT' \sub \tT \in A^H_B \} = m^\iB_{\iB'} H_\iB(f)^*$,
so $\sgen{\gG^\sharp[O]}$ is generated by vectors 
$v^{Hf} \in (\mb{Z}^Q)^O$
where $H \in \mc{H}$, $f \sub [q]$ with $i_{\mc{P}}(f)=\iB'$
and for each $\psi \in O$, $B \in Q$ we have 
$(v^{Hf}_{\psi})_B = m^\iB_{\iB'} H_\iB(f)^*$, 
where $\iB=i_{\mc{P}}(B)$.
Thus all vectors in $\sgen{\gG^\sharp[O]}$ are obtained
from vectors $H(f)^*$ with $H \in \mc{H}$ and 
$i_{\mc{P}}(f) = i_{\mc{P}'}(e)$ by the same transformation
that maps $G(e)^*$ to $((G^*)^\sharp)^O$.
As $G$ is $(\mc{H},\mc{P})$-divisible we deduce
$((G^*)^\sharp)^O \in \sgen{\gG^\sharp[O]}$ 
for any $O \in \Phi/\Ss$, as required. \qed

\section{Directed hypergraphs} \label{sec:di}

Our second illustration of Theorem \ref{main}
will be to decompositions of directed hypergraphs.

\begin{defn} 
Let $R$ be a set.
An $R$-graph on $V$ is a set $G$ of injections from $R$ to $V$.
We call the elements of $G$ arcs.
If $R=[r]$ we call $G$ an $r$-digraph.
We say $G$ is simple if $(Im(e): e \in G)$ are all distinct.
A copy of an $R$-graph $H$ in an $R$-graph $G$
is defined by an injection $\phi:V(H) \to V(G)$
such that $\phi H := \{ \phi \circ e: e \in H \} \sub G$.
An $H$-decomposition of $G$ is a partition of $G$
into copies of $H$.
\end{defn} 

Note that if $r=2$ then a $2$-digraph is equivalent
to a digraph in the usual sense: we can think of an injection 
$f:[2] \to V$ as an arc directed from $f(1)$ to $f(2)$.

We will restrict our attention to $H$-decomposition problems
in which $H$ is simple; otherwise we obtain a non-elementary 
functional decomposition problem, which has arithmetic structure,
and to which Theorem \ref{main} does not apply.

Next we will state an example of our later theorem
on $r$-digraph decompositions. Let $KD^r_n$ denote
the complete $r$-digraph on $[n]$, i.e.\ each of the
$(n)_r = r!\tbinom{n}{r}$ injections from $[r]$ to $[n]$ is an arc.
The $r$-digraph tight $q$-cycle $\car^r_q$ has vertex set $[q]$
and arc set $\{ \phi_j: j \in [q] \}$ with each $\phi_j(i)=i+j$,
where addition wraps (we identify $q+i$ with $i$).

\begin{theo} \label{KDC}
Suppose $q>r \ge 2$ and $n>n_0(q)$ with $q \mid (n)_r$.
Then $KD^r_n$ has a $\car^r_q$-decomposition.
\end{theo}

Now we will describe the divisibility conditions
in the general setting, and then illustrate them
in the case $H=\car^r_q$.

\begin{defn} \label{didiv}
Let $G$ be an $r$-digraph on $[n]$
and $H$ be an $r$-digraph on $[q]$.

Given an injection $f:R' \to [n]$ with $R' \sub R$,
we let $G\mid_f = \{e \in G: e\mid_{R'}=f \}$.
The neighbourhood of $f$ in $G$ is the $(R \sm R')$-graph
$G(f) = \{ e\mid_{R \sm R'}: e \in G\mid_f\}$.
The degree of $f$ in $G$ is $|G(f)|$.

We write $I^s_t$ for the set of 
injections $\pi: [s] \to [t]$.
For $\psi \in I^i_n$ we define the degree vector 
$G(\psi)^* \in \mb{N}^{I^i_r}$
by $G(\psi)^*_\pi = |G(\psi \pi^{-1})|$.

We say $G$ is $H$-divisible if 
$G(\psi)^* \in \sgen{ H(\tT)^* : \tT \in I^i_q }$
for all $0 \le i \le r$, $\psi \in I^i_n$.
\end{defn}

Now we illustrate Definition \ref{didiv} in the case $H=\car^r_q$. 
For example, suppose $r=2$, so $H$ and $G$ are digraphs.
Writing $\es$ for the element of $I^0_n$,
we have $G(\es)^* = (|G|)$ and $H(\es)^*=(|H|)=(q)$,
so the $0$-divisibility condition is $q \mid |G|$.
Next, for $\psi \in I^1_n$, writing $x=\psi(1) \in [n]$,
we have $G(\psi)^* = (d^+_G(x),d^-_G(x))$, where 
$d^+_G(x)=|G(\psi)|$ is the number of arcs with $1 \mapsto x$
and $d^-_G(x)=|G(\psi \circ (1 \mapsto 2)^{-1})|$ 
is the number of arcs with $2 \mapsto x$.
Also, for $\tT \in I^1_q$, writing $a=\tT(1) \in [q]$,
we have $H(\tT)^*=(d^+_H(a),d^-_H(a))=(1,1)$, so the 
$1$-divisibility condition is that $G$ is vertex-regular,
i.e.\ $d^+_G(x)=d^-_G(x)$ for all $x \in [n]$.
Finally, for $\psi \in I^2_n$, $\tT \in I^2_q$
writing $x_i=\psi(i)$, $a_i=\tT(i)$, we have 
$G(\psi)^* = (1_{x_1x_2 \in G}, 1_{x_2x_1 \in G})$
and $H(\tT)^* = (1_{a_1a_2 \in H}, 1_{a_2a_1 \in H})$,
so the $2$-divisibility condition holds trivially.
Next we describe the general $\car^r_q$-divisibility conditions 
(proved in Lemma \ref{Cdiv} below). 

\begin{defn}
We define an equivalence relation $\sim$ on each $I^i_r$
with $i \le r$ by $\tT \sim \tT'$ if for some $c \in \mb{Z}$ 
we have $\tT'(j)=\tT(j)+c$ for all $j \in [i]$
(where addition does not wrap).
We say that $G$ is shift regular if 
$G(\psi)^*_\tT = G(\psi)^*_{\tT'}$ 
whenever $\psi \in I^i_n$ and $\tT \sim \tT'$.
\end{defn}

We note that $G =  KD^r_n$ is shift regular, indeed
$G(\psi)^*_\tT = (n)_r/(n)_i$ for any $\tT \in I^i_r$, $\psi \in I^i_n$.
We also note that there is redundancy (symmetry)
in the above definitions. Indeed, for $\psi \in I^i_n$,
$\sS \in S_i$, $\pi \in I^i_r$ we have $G(\psi \sS)^*_\pi
= |G(\psi \sS \pi^{-1})| = G(\psi)^*_{\pi \sS^{-1}}$,
i.e.\ $G(\psi \sS)^* = G(\psi)^* \sS$,
where $S_i$ acts on $I^i_n$ by 
$\psi \mapsto \psi \sS = \psi \circ \sS$
and on $\mb{N}^{I^i_r}$ by $(v \sS)_\pi = v_{\pi \sS^{-1}}$.
Note that the latter is a right action as
$(v (\sS\tau))_\pi = v_{\pi (\sS\tau)^{-1}}
= v_{\pi \tau^{-1} \sS^{-1} } = (v \sS)_{\pi \tau^{-1}}
= ((v \sS)\tau)_\pi$.
For any expression $G(\psi)^* = \sum_\tT n_\tT H(\tT)^*$
with $n \in \mb{Z}^{I^i_q }$ we have
$G(\psi \sS)^* = G(\psi)^* \sS = \sum_\tT n_\tT H(\tT)^* \sS
= \sum_\tT n_\tT H(\tT \sS)^*$, so it suffices to check
$H$-divisibility on a system of coset 
representatives for the action of $S_i$ on $I^i_n$.
Furthermore, as $\tT \sim \tT'$ iff $\tT \sS \sim \tT' \sS$,
and as $G(\psi)^*_{\tT\sS} = |G(\psi (\tT \sS)^{-1})|
= G(\psi \sS^{-1})^*_{\tT} $, 
it suffices to check shift regularity on a 
system of coset representatives for the action of $S_i$ 
on $I^i_q$, e.g.\ all order-preserving elements.

\begin{lemma} \label{Cdiv}
$G$ is $\car^r_q$-divisible 
iff $G$ is shift regular and $q \mid |G|$.
\end{lemma}

\nib{Proof.}
The $0$-divisibility condition is $q \mid |G|$. Fix $0<i \le r$.
We classify the degree vectors $H(\tT)^*$ with $\tT \in I^i_q$. 
Note that $H(\tT)^*$ is the all-0 vector 
unless $Im(\tT)$ is contained in a cyclic interval of length $r$.
By the cyclic symmetry of $\car^r_q$ we have $H(\tT)^*=H(\tT+c)^*$
for any $c \in [q]$, defining $\tT+c \in I^i_q$ by 
$\tT(j)=\tT'(j)+c$ (where addition wraps).
Thus we can assume $R := Im(\tT) \sub [r]$, i.e.\ $\tT \in I^i_r$.
Note that $id_{[r]}$ is the unique arc of $H$ containing
$id_R$, so $1=|H(id_R)|=H(\tT)^*_\tT$. 
Similarly, for each $c \in \mb{Z}$ such that $R+c \sub [r]$
(where addition does not wrap), $id_{[r]}-c$ is the unique arc of $H$ 
containing $id_{R+c}-c$, so $1=|H(id_{R+c}-c)|=H(\tT)^*_{\tT+c}$.
All other coordinates of $H(\tT)^*$ are zero.
We deduce that $H(\tT)^*=H(\tT')^*$ if $\tT \sim \tT'$,
or otherwise $H(\tT)$ and $H(\tT')^*$ have disjoint support.
Thus $G(\psi)^* \in \sgen{ H(\tT)^* : \tT \in I^i_q }$
iff $G$ is constant on the support of each $H(\tT)^*$,
i.e.\ $G$ is shift regular. \qed

\medskip

Given Lemma \ref{Cdiv}, the case $H=\car^r_q$ of the following result
implies Theorem \ref{KDC}.

\begin{theo} \label{KD}
Suppose $H$ is a simple $r$-digraph on $[q]$ and $n>n_0(q)$ is large.
Then $KD^r_n$ has an $H$-decomposition iff it is $H$-divisible.
\end{theo}

We will deduce Theorem \ref{KD} from a more general result
in which we replace $KD^r_n$ by any $r$-digraph $G$ 
supported in a $[q]$-complex $\Phi$ that satisfies
certain extendability and regularity conditions.
The regularity condition is similar to those
used earlier in the paper:

\begin{defn} \label{def:direg}
Let $\Phi$ be a $[q]$-complex on $[n]$, 
$H$ be an $r$-digraph on $[q]$
and $G$ be an $r$-digraph on $[n]$.
We say $G$ is $(H,c,\oO)$-regular in $\Phi$ if there are
$y_\phi \in [\oO n^{r-q},\oO^{-1} n^{r-q}]$ 
for each $\phi \in \Phi_q$ with $\phi H \sub G$
so that $\sum_\phi y_\phi \phi H = (1 \pm c)G$.
\end{defn}

Next we introduce some notation for the extendability 
condition and illustrate it for digraphs.

\begin{defn} \label{def:diext}
With notation as in Definition \ref{def:direg},
let $Q^H$ be the set of $B \in Q = \tbinom{[q]}{r}$ 
such that there is some $\tT_B \in H$ with $Im(\tT_B)=B$.
Suppose $H$ is simple, so that each $\tT_B$ is unique.
Define $G^H \sub \Phi_r$ by $G^H_B = 
\{ \psi \circ \tT_B^{-1}: \psi \in G \}$ 
if $B \in Q^H$ or $G^H_B = \Phi_B$ otherwise. 
\end{defn}

\nib{Examples.}
Let $q=3$, $r=2$, $G$ be a digraph on $[n]$ and
$\Phi$ be the complete $[3]$-complex on $[n]$.
\begin{enumerate}
\item Let $H = \{ (1 \mapsto 1, 2 \mapsto 2),
(1 \mapsto 2, 2 \mapsto 3), (1 \mapsto 3, 2 \mapsto 1) \}$
be a cyclic triangle. For each $i \in [3]$ we have
$G^H_{\{i,i+1\}} = \{ (i \mapsto x, i+1 \mapsto y):
xy=(1 \mapsto x, 2 \mapsto y) \in G \}$
(interpreting $i+1$ mod $3$).

If $(\Phi,G^H)$ is $(\oO,h)$-extendable
then for any disjoint sets $S_i \sub T_i$, $i \in [3]$
of size at most $h$ and injection 
$\phi: S:=\bigcup_{i=1}^3 S_i \to [n]$
there are at least $\oO n^{|T \sm S|}$
injections $\phi^+: T:=\bigcup_{i=1}^3 T_i \to [n]$
extending $\phi$ such that for any $i \in [3]$,
$x_i \in T_i$, $x_{i+1} \in T_{i+1}$ (addition mod $3$) 
with $x_ix_{i+1} \not\sub S$ we have 
$(i \mapsto \phi^+(x_i), i+1 \mapsto \phi^+(x_{i+1}))
\in G^H_{\{i,i+1\}}$, i.e.\ $\phi^+(x_i)\phi^+(x_{i+1}) \in G$.

This is roughly equivalent to the following property:
say that $G$ is fully $(\oO,h)$-extendable
if for any disjoint $A,B \sub [n]$ of size at most $h$
there are at least $\oO n$ vertices $c$
such that $ca \in G$ for all $a \in A$
and $bc \in G$ for all $b \in B$.
Indeed, if $(\Phi,G^H)$ is $(\oO,h)$-extendable
then $G$ is fully $(\oO,h)$-extendable
(take $S_1=T_1=A$, $S_2=T_2=B$, $S_3=\es$, $|T_3|=1$),
and conversely, if $G$ is fully $(\oO,h)$-extendable
then $(\Phi,G^H)$ is $(\oO^{3h},h)$-extendable
(construct $\phi^+$ one vertex at a time).

\medskip

\item Now let $H = \{ (1 \mapsto 1, 2 \mapsto 2),
(1 \mapsto 1, 2 \mapsto 3) \}$ be an outstar 
of degree two. For $i=2,3$ we have
$G^H_{1i} = \{ (1 \mapsto x, i \mapsto y):
xy \in G \}$, and $G^H_{23}=\Phi_{23}$ is complete.
If $(\Phi,G^H)$ is $(\oO,h)$-extendable
then given $S_i$, $T_i$ and $\phi$ as above,
there are at least $\oO n^{|T \sm S|}$
extensions $\phi^+$ such that for any $i=2,3$,
$x_1 \in T_1$, $x_i \in T_i$ with $x_1x_i \not\sub S$ 
we have $\phi^+(x_1)\phi^+(x_i) \in G$.

This is roughly equivalent to the following property:
say that $G$ is directedly $(\oO,h)$-extendable
if for any $A \sub [n]$ of size at most $h$
there are at least $\oO n$ vertices $c$
such that $ca \in G$ for all $a \in A$,
and at least $\oO n$ vertices $c$
such that $ac \in G$ for all $a \in A$.
\end{enumerate}

\medskip

The rough equivalence illustrated in the previous examples 
takes the following general form:
if $(\Phi,G^H)$ is $(\oO,h)$-extendable
then $(\Phi,G)$ is $(\oO,h,H)$-vertex-extendable
(as in the next definition),
and conversely, if $G$ is $(\oO,h,H)$-vertex-extendable
then $(\Phi,G^H)$ is $(\oO^{qh},h)$-extendable.

\begin{defn} \label{def:diext=}
With notation as in Definition \ref{def:direg},
we say $(\Phi,G)$ is $(\oO,h,H)$-vertex-extendable
if for any $x \in [q]$ and disjoint sets $A_i$, 
$i \in [q] \sm \{x\}$ of size at most $h$ such that
$(i \mapsto v_i: i \in [q]\sm \{x\}) \in \Phi$
whenever each $v_i \in A_i$, there are at least $\oO n$ 
vertices $v \in \Phi^\circ_x$ such that
\begin{enumerate}
\item $(i \mapsto v_i: i \in [q]) \in \Phi$
whenever $v_x=v$ and $v_i \in A_i$ for each $i \ne x$,
\item for each arc $\tT$ of $H$ with $x \in Im(\tT)$,
we have all arcs $(i \mapsto v_i: i \in [r])$ in $G$
where $v_j=v$ for $j=\tT^{-1}(x)$
and $v_i \in A_{\tT(i)}$ for all $i \ne j$. 
\end{enumerate}
\end{defn}

The following theorem when $\Phi$ and $G$ 
are complete implies Theorem \ref{KD}. 
Indeed, extendability is clear, and for regularity we let
$y_\phi = |H|^{-1} (n)_r/(n)_q$ for each $\phi \in I^q_n$,
so that for each $\psi \in I^r_n$ we have
$\sum_\phi y_\phi (\phi H)_\psi 
= \sum_{\tT \in H}  |H|^{-1} (n)_r (n)_q^{-1} 
|\{\phi: \psi=\phi\tT\}| = 1$.

\begin{theo} \label{diext}
Let $H$ be a simple $r$-digraph on $[q]$,
$G$ be an $r$-digraph on $[n]$ and $\Phi$ be an 
$(\oO,h)$-extendable $S_q$-adapted $[q]$-complex on $[n]$
where $n>n_0(q)$ is large, 
$n^{-\dD}<\oO<\oO_0(q)$ is small and $c=\oO^{h^{20}}$.
Suppose $G$ is $(H,c,\oO)$-regular in $\Phi$
and $(\Phi,G^H)$ is $(\oO,h)$-extendable.
Then $G$ has an $H$-decomposition in $\Phi_q$
iff $G$ is $H$-divisible.
\end{theo}

To deduce this from Theorem \ref{main} we will use 
the following equivalent encoding.

\begin{defn} \label{def:di*}
Given an injection $f:[r] \to X$, we write 
$f^{r \to q}$ for the set of all $f \circ \pi^{-1}$
where $\pi:[r] \to [q]$ is order-preserving.
Given an $r$-digraph $G$, we let $G^{r \to q}$ be the
(disjoint) union of all $f^{r \to q}$ with $f \in G$.
\end{defn}

\begin{lemma} \label{equivdi*}
Let $H$ and $G$ be $r$-digraphs,
$H^* = H^{r \to q}$ and $G^* = G^{r \to q}$.
Then an (integral) $H$-decomposition of $G$ is equivalent to
an (integral) $H^*$-decomposition of $G^*$. 
\end{lemma}

\nib{Proof.}
We associate any $H$-decomposition $\mc{H}$ of $G$
with an $H^*$-decomposition $\mc{H}^*$ of $G^*$, associating 
each $\phi H \in \mc{H}$ with $\phi H^* \in \mc{H}^*$.
Then $e \in \phi H$ iff $e^{r \to q} \sub \phi H^*$,
as if $e = \phi\tT$ for some $\tT \in H$
and $e\pi^{-1} \in e^{r \to q}$ then $e\pi^{-1} = \phi \tT^*$, 
where $\tT^* = \tT \pi^{-1} \in H^*$, and conversely.
The same proof applies to integral decompositions.  \qed

\medskip

\nib{Proof of Theorem \ref{diext}.}
Let $H^* = H^{r \to q}$ and $G^* = G^{r \to q}$.
Let $\mc{A}=\{A\}$ with $A=S_q^\le$ and
$\gG \in \mb{Z}^{A_r}$ where each $\gG_\tT = 1_{\tT \in H^*}$.
Then a $\gG(\Phi)$-decomposition of $G^*$
is equivalent to an $H^*$-decomposition of $G^*$, and so
(by Lemma \ref{equivdi*}) to an $H$-decomposition of $G$.

Next we claim that $\gG$ is elementary. To see this, 
we describe the type vectors $\gG^\tT \in \{0,1\}^{(S_q)^B}$ 
for $\tT \in A_B$, $B \in Q = \tbinom{[q]}{r}$.
If $\gG^\tT \ne 0$ then we can write $\tT = \tT_0 \pi_0 \sS_0$
with $\tT_0 \in H$, $\pi_0 \in S_r$ and $\sS_0 \in Bij(B,[r])$ 
order-preserving; this expression is unique, 
as $\tT_0$ is determined by $\tT$ (as $H$ is simple).
For any $\sS \in \Ss^B$ we have $\gG^\tT_\sS = \gG_{\tT\sS}$ 
equal to $1$ iff $\sS = \sS_0^{-1} \pi_0^{-1} \pi^{-1}$ 
where $\pi:[r] \to [q]$ is order-preserving.
Thus there are $r!+1$ types: the $0$ type, and types $t^{\pi_0}$ 
for each $\pi_0 \in S_r$, describing the $r!$ possible arcs 
with any given image. The supports of the $t^{\pi_0}$ are 
mutually disjoint, so $\gG$ is elementary, as claimed.

The atom decomposition is 
$G^* = \sum_{e \in G} e^*$, where $e^* = e^{r \to q}$.
As $G$ is $(H,c,\oO)$-regular in $\Phi$, we have
$\sum_\phi y_\phi \phi H = (1 \pm c) G$ 
(equivalently, $\sum_\phi y_\phi \phi H^* = (1 \pm c) G^*$)
for some $y_\phi \in [\oO n^{r-q},\oO^{-1} n^{r-q}]$ 
for each $\phi \in \Phi_q$ with $\phi H \sub G$
(equivalently, $\phi H^* \sub G^*$).
For any such $\phi$ we have $\gG(\phi) \le_\gG G^*$,
so $\phi \in \mc{A}(\Phi,G^*)$. Also, for any
$B \in Q$, $\psi \in \Phi_B$ and $0 \ne t \in T_B$,
say with $t$ supported on the set of all $\tau^{-1} \pi^{-1}$
where $\pi:[r] \to [q]$ is order-preserving,
we have $\pl^t y_\psi 
= \sum \{ y_\phi : t_\phi(\psi)=t \}
= \sum \{ y_\phi: \psi\tau \in \phi H^* \}
= (1 \pm c) G^*_{\psi\tau}
= (1 \pm c)(G^*)^t_\psi$,
so $G^*$ is $(\gG,c,\oO)$-regular.

Next we consider extendability. We have
$\gG[G^*] = \{ \psi \in \Phi_r : \gG(\psi) \le_\gG G^* \}$,
so $\psi \in \Phi_B$ is in $\gG[G^*]$ iff 
(a) no arc in $H$ has image $B$, or
(b) $\psi \tT_B \in G$ for
the unique arc $\tT_B$ in $H$ with $Im(\tT_B)=B$.
Let $E=(J,F,\phi)$ be any $\Phi$-extension of rank $h$
and $J' \sub J_r \sm J[F]$. 
As $(\Phi,G^H)$ is $(\oO,h)$-extendable
we have $X_{E,J'}(\Phi,G^H) > \oO n^{v_E}$.
Consider any $\phi^+ \in X_{E,J'}(\Phi,G^H)$.
For any $\psi \in J'_B$ we have 
$\phi^+\psi  \in G^H_B$, 
so $\phi^+\psi \tT_B  \in G$, 
so $\phi^+\psi \in \gG[G^*]$.
Thus $\phi^+ \in X_{E,J'}(\Phi,\gG[G^*])$,
so $(\Phi,\gG[G^*])$ is $(\oO,h)$-extendable.

To deduce the theorem from Theorem \ref{main},
it remains to consider divisibility.
By Lemma \ref{lattice} we have
$\bgen{\gG(\Phi)} = \mc{L}_\gG(\Phi)$.
By Definition \ref{def:L} we need to show
that  $G$ is $H$-divisible iff
$((G^*)^\sharp)^O \in \sgen{\gG^\sharp[O]}$
for any orbit $O \in \Phi/S_q$. 
To describe $((G^*)^\sharp)^O \in (\mb{Z}^Q)^O$, recall that
if $\psi' \in O \cap \Phi_{B'}$ then
$((G^*)^\sharp_{\psi'})_B)$  is the number of 
$\psi \in G^* \cap \Phi_B$ with $\psi\mid_{B'}=\psi'$.
We can assume $B' \sub B$, otherwise this number is $0$.
Let $\pi_B:[r] \to B$ be order-preserving and $R=\pi_B^{-1}(B')$.
Then $\psi \in G^* \cap \Phi_B$ iff $\psi \pi_B \in G$,
and $\psi\mid_{B'}=\psi'$ iff $(\psi \pi_B)\mid_R=\psi' \pi_B$,
so $((G^*)^\sharp_{\psi'})_B=|G(\psi' \pi_B)|$.
Similarly, to describe $\sgen{\gG^\sharp[O]}$, recall that
it is generated by vectors $\gG^\sharp(\phi) \in (\mb{Z}^Q)^O$
where if $\psi' = \phi \tT'$ with $\tT' \in A_{B'}$ then
$(\gG^\sharp(\phi)_{\psi'})_B = (\gG^\sharp_{\tT'})_B$
is the number of $\tT \in H^*_B$ with $\tT\mid_{B'}=\tT'$,
which is $|H(\tT' \pi_B)|$.

Now fix $\psi \in O \cap \Phi_{[i]}$, where $O \in \Phi_i/S_q$.
As $G$ is $H$-divisible, there is $n \in \mb{Z}^{I^i_q}$ with
$G(\psi)^* = \sum_{\tT} n_{\tT} H(\tT)^*$.
Writing $\phi = \psi \tT^{-1}$, we claim that 
$((G^*)^\sharp)^O = \sum_{\tT} n_{\tT} \gG^\sharp(\phi)$.
To see this, note that it suffices to prove
$((G^*)^\sharp)^O_{[r]} = \sum_{\tT} n_{\tT} \gG^\sharp(\phi)_{[r]}$,
as $((G^*)^\sharp_{\psi'})_B = |G(\psi' \pi_B)|
= ((G^*)^\sharp_{\psi'\pi_B})_{[r]}$ and
$(\gG^\sharp(\phi)_{\psi'})_B 
= (\gG^\sharp_{\phi^{-1}\psi'})_B
= |H(\phi^{-1}\psi'\pi_B)|
= (\gG^\sharp(\phi)_{\psi'\pi_B})_{[r]}$.
Now for any $\psi' \in O \cap \Phi_R$ with $R \sub [r]$,
writing $\pi = (\psi')^{-1} \psi \in I^i_r$, we have
$((G^*)^\sharp_{\psi'})_{[r]} 
= |G(\psi')| = G(\psi)^*_\pi
= \sum n_{\tT} H(\tT)^*_\pi$,
where each $H(\tT)^*_\pi = |H(\tT \pi^{-1})|
= (\gG^\sharp_{\tT \pi^{-1}})_{[r]}
= (\gG^\sharp(\phi)_{\psi'})_{[r]}$,
so $((G^*)^\sharp_{\psi'})_{[r]} 
= \sum n_{\tT} (\gG^\sharp(\phi)_{\psi'})_{[r]}$. \qed

\section{All of the above}\label{sec:all}

For use in future applications (e.g.\ \cite{KS}),
in this section we present a general theorem that simultaneously 
allows for the various flavours of decomposition considered
in this paper (generalised partitions, colours and directions).
We start with a definition that generalises our previous setting
of simple $r$-digraphs to allow for colours, 
index vectors with respect to a partition, and different
types of `generalised arcs'; it is followed by some 
illustrative examples.

\begin{defn} \label{def:canonical}
Let $\mc{P}=(P_1,\dots,P_t)$ be a partition of $[q]$ such that 
if $x \in P_j$, $x' \in P_{j'}$, $j<j'$ then $x<x'$.
Let $\mc{H}$ be a family of $[D]$-edge-coloured 
$r$-digraphs on $[q]$. For $\iB \in \mb{N}^t$
with $\sum_{j=1}^t i_j=r$ and $j \in [t]$ we define 
a partition $R(\iB)=(R(\iB)_1,\dots,R(\iB)_t)$ of $[r]$ 
so that each $|R(\iB)_j|=i_j$ and $x<x'$ 
whenever $x \in R(\iB)_j$, $x' \in R(\iB)_{j'}$, $j<j'$.  
Suppose there are vectors $\iB^d \in \mb{N}^t$
and permutation groups $\LL^d_j$ on $R(\iB)_j$
for all $d \in [D]$ and $j \in [t]$
such that if $H \in \mc{H}$ and $\tT \in H^d$ then 
\begin{enumerate}
\item each $\tT(R(\iB)_j) \sub P_j$
(so\footnote{Recall index vectors from Definition \ref{def:HPblowup}.} $i_{\mc{P}}(Im(\tT))=\iB^d$), and
\item
for $\tT' \in Bij([r],Im(\tT))$ we have
$\tT' \notin H \sm H^d$, and $\tT' \in H^d$ 
iff $\tT^{-1} \tT' \in \LL^d := \prod_j \LL^d_j$. 
\end{enumerate}
We say that $\mc{H}$ is $(\mc{P},\LL)$-canonical,
where $\LL := (\LL^d: d \in D)$.
\end{defn}

\nib{Examples.}
\begin{enumerate}
\item Let $q=3$, $r=2$ and $t=1$, 
so $\mc{P}=([3])$ and $R(2)=([2])$.
Let $D=2$ and $\mc{H}=\{H\}$, where 
$H^1 = \{ (1 \mapsto 1, 2 \mapsto 2),
(1 \mapsto 1, 2 \mapsto 3) \}$ and
$H^2 = \{ (1 \mapsto 2, 2 \mapsto 3),
(1 \mapsto 3, 2 \mapsto 2) \}$.
Then $\mc{H}$ is canonical 
with $\LL^1_1 = \{id\}$ and 
$\LL^2_1 = S_2 = \{id,(12)\}$.
One can interpret $H$ as a mixed triangle,
with arcs from $1$ to $2$ and $1$ to $3$
and an undirected edge between $2$ and $3$.
In this interpretation, we are free to ignore 
the colours, as they do not affect whether 
a mixed graph $G$ has an $H$-decomposition
(the role of the colours is to ensure 
that under the encoding by arcs,
an undirected edge encoded by two arcs
cannot be decomposed into two actual arcs). 
In general, we think of an atom in some colour $d$
as a `generalised arc', which is encoded by some set
of arcs invariant under the action of $\LL^d$ on $[r]$. 
An actual arc corresponds to the
case $\LL^d=\{id\}$ and an undirected edge
to the case that each $\LL^d_j=Sym(R(\iB)_j)$.

\medskip

\item Let $q=3$, $r=2$, $t=1$, $D=2$
and $\mc{H}=\{H\}$, where 
$H^1 = \{ (1 \mapsto 1, 2 \mapsto 2),
(1 \mapsto 2, 2 \mapsto 3) \}$ and
$H^2 = \{ (1 \mapsto 3, 2 \mapsto 1)\}$.
Then $\mc{H}$ is canonical 
with $\LL^1 = \LL^2 = \{id\}$.
One can interpret $H$ as a two-coloured cyclic
directed triangle, with arcs of colour $1$ 
from $1$ to $2$ and $2$ to $3$, and an
arc of colour $2$ from $3$ to $1$.

\medskip

\item Let $q=3$, $r=2$, $t=2$,
$\mc{P}=(\{1,2\},\{3\})$, $D=3$
and $\mc{H}=\{H\}$, where 
$H^1 = \{ (1 \mapsto 1, 2 \mapsto 2) \}$,
$H^2 = \{ (1 \mapsto 1, 2 \mapsto 3) \}$ and
$H^3 = \{ (1 \mapsto 2, 2 \mapsto 3) \}$.
We have $\iB^1=(2,0)$, $R((2,0))=([2],\es)$,
$\iB^2=\iB^3=(1,1)$, $R((1,1))=(\{1\},\{2\})$
and $\LL^1=\LL^2=\LL^3=\{id\}$.
One possible uncoloured interpretation of $H$ is
as a cyclic triangle $1 \to 2 \to 3 \to 1$
under the vertex partition $\mc{P}$.
Here we are taking the natural interpretation 
of the colour $3$ arc from $2$ to $3$
and the opposite interpretation 
of the colour $2$ arc from $1$ to $3$,
instead thinking of it as an arc from $3$ to $1$.
Changing the direction of all arcs of colour $2$
in both $H$ and $G$ has no effect on 
whether $G$ has an $H$-decomposition,
so this interpretation is equivalent to the
natural interpretation in which we retain
the given colours and directions.
This illustrates the fact that in general there is 
no loss of generality from the assumption that the partitions
$\mc{P}$ and $R(\iB)$ respect the orders of $[q]$ 
and $[r]$, as we are free to interpret different 
colours as encoding arcs with alternative partitions.
We also note that there is no loss of generality
in assuming that the index of an edge is determined
by its colour (and indeed, we could have done so
earlier in the paper).
\end{enumerate}

\medskip

For the main result of this section
we adopt the setting of the following definition
(see below for how it applies to the above examples).

\begin{defn} \label{def:master}
Let $\mc{H}$ be a $(\mc{P},\LL)$-canonical family 
of $[D]$-edge-coloured $r$-digraphs on $[q]$.
We identify each $H \in \mc{H}$ with a vector
$H \in (\mb{N}^D)^{I^r_q}$,
where each $(H_f)_d = 1_{f \in H^d}$.

Let $\Ss$ be the group of all $\sS \in S_q$ 
with all $\sS(P_i)=P_i$.
Let $\Phi$ be an exactly $\Ss$-adapted $[q]$-complex
with $V(\Phi)=[n]$ and parts $\mc{P}'=(P'_1,\dots,P'_t)$, where each 
$P'_i = \{ \psi(j): j \in P_i,\ \psi \in \Phi_{\{j\}} \}$.
For $\phi \in \Phi_q$ and $H \in \mc{H}$ we define 
$\phi H \in (\mb{N}^D)^{\Phi_{[r]}}$
by $(\phi H)_{\phi f} = H_f$.
Let $\mc{H}(\Phi) = \{ \phi H:
 \phi \in \Phi_q, H \in \mc{H} \}$.

Let $G \in (\mb{N}^D)^{\Phi_{[r]}}$ be an $r$-multidigraph 
$[D]$-edge-coloured as $G = \cup_{d \in [D]} G^d$.

We call $\mc{H'} \sub \mc{H}(\Phi)$ with 
$\sum \mc{H'} = G$ an $H$-decomposition of $G$ in $\Phi$.

We call $\Psi \in \mb{Z}^{\mc{H}(\Phi)}$ with
$\sum_{H'} \Psi_{H'} H' = G$
an integral $H$-decomposition of $G$ in $\Phi$.

For $\psi \in I^i_n$ (injections $[i] \to [n]$)
and $\tT \in I^i_q$ write
$i_{\mc{P}'}(\psi)=i_{\mc{P}'}(Im(\psi))$
and $i_{\mc{P}}(\tT)=i_{\mc{P}}(Im(\tT))$.

For $\psi \in I^i_n$ we define the degree vector 
$G(\psi)^* \in \mb{N}^{[D] \times I^i_r}$
by $G(\psi)^*_{d\pi} = |G^d(\psi \pi^{-1})|$.

Similarly, for $\tT \in I^i_q$ we define
$H(\tT)^* \in \mb{N}^{[D] \times I^i_r}$
by $H(\tT)^*_{d\pi} = |H^d(\tT \pi^{-1})|$.
For $\iB' \in \mb{N}^t$ we let 
$H\sgen{\iB'} = \sgen{H(\tT)^*: i_{\mc{P}}(\tT)=\iB'}$.
We say $G$ is $\mc{H}$-divisible (in $\Phi$) if
$G(\psi)^* \in H\sgen{\iB'}$ whenever $i_{\mc{P}'}(\psi)=\iB'$.

We say $G$ is $(\mc{H},c,\oO)$-regular in $\Phi$ if there are
$y^H_\phi \in [\oO n^{r-q},\oO^{-1} n^{r-q}]$ for each 
$H \in \mc{H}$, $\phi \in \Phi_q$ with $\phi H \le G$
so that $\sum \{ y^H_\phi \phi H \} = (1 \pm c)G$.

For each $H \in \mc{H}$ and $B \in Q$ 
fix any $\tT_B \in H$ with $Im(\tT_B)=B$ if one exists.

For each $d \in [D]$ let $(G^H)^d = \bigcup 
\{ \psi \circ \tT_B^{-1}: \tT_B \in H^d,\ G^d_\psi>0 \}$. 

We say that $(\Phi,G^H)$ is $(\oO,h)$-extendable
if $(\Phi,((G^H)^d: d \in [D]))$ is $(\oO,h)$-extendable.
\end{defn}

\nib{Examples.}
\begin{enumerate}
\item Recall the example of the mixed triangle:
$q=3$, $r=2$, $t=1$, $D=2$, $\mc{H}=\{H\}$, 
$H^1 = \{ (1 \mapsto 1, 2 \mapsto 2),
(1 \mapsto 1, 2 \mapsto 3) \}$,
$H^2 = \{ (1 \mapsto 2, 2 \mapsto 3),
(1 \mapsto 3, 2 \mapsto 2) \}$,
$\LL^1_1 = \{id\}$, $\LL^2_1 = \{id,(12)\}$.
Let $\Phi$ be the complete $[3]$-complex on $[n]$
and $G \in (\mb{N}^2)^{\Phi_2}$ be 
a $[2]$-edge-coloured $2$-multidigraph.
For the $2$-divisibility condition
we consider any $\psi \in I^2_n$, so that
$G(\psi)^* \in \mb{N}^{[2] \times I^2_2}$.
Ordering coordinates as
$(1,id)$, $(1,(12))$, $(2,id)$, $(2,(12))$ we have 
$G(\psi)^* = ( G^1_\psi, G^1_{\psi \circ (12)},
 G^2_\psi, G^2_{\psi \circ (12)} )$.
The possible $H(\tT)^*$ with $\tT \in I^2_3$
are $(1,0,0,0)$, $(0,1,0,0)$ and $(0,0,1,1)$.
Thus the $2$-divisibility condition is that
$G^2_\psi = G^2_{\psi \circ (12)}$ 
for all $\psi \in I^2_n$, i.e.\ arcs of colour $2$
always come in opposite pairs (which we interpret
as an edge when we think of $G$ as a mixed multigraph).
As for $0$-divisibility,
writing $\es$ for the function with empty domain,
we have $G(\es)^* = (|G^1|,|G^2|)$
and $H(\es)^*=(|H^1|,|H^2|)=(2,2)$,
so we need $|G^1|=|G^2|$. In terms of mixed multigraphs,
we need twice as many arcs as edges 
(each edge corresponds to a pair of arcs in $G^2$).

For the $1$-divisibility conditions,
consider any $\psi \in I^1_n$, 
so $G(\psi)^* \in \mb{N}^{[2] \times I^1_2}$.
Let $x=Im(\psi) \in [n]$. Ordering coordinates as
$(1,1 \mapsto 1)$, $(1,1 \mapsto 2)$, 
$(2,1 \mapsto 1)$, $(2,1 \mapsto 2)$ we have 
$G(\psi)^* = (|G^1(1 \mapsto x)|, |G^1(2 \mapsto x)|,
|G^2(1 \mapsto x)|, |G^2(2 \mapsto x)| )
= ( d_G^+(x), d_G^-(x), d_G(x), d_G(x) )$,
where in the mixed graph interpretation
$d_G^\pm(x)$ denote in/outdegrees in arcs
and $d_G(x)$ denotes degree in edges.
We have $H(1 \mapsto 1)^* = (2,0,0,0)$ and
$H(1 \mapsto 2)^* = H(1 \mapsto 3)^* = (0,1,1,1)$.
Thus the $1$-divisibility conditions are that 
each outdegree $d_G^+(x)$ is even
and each $d_G(x)=d_G^-(x)$.

Now consider extendability. We have $(G^H)^1_{12} = G^1$,
$(G^H)^1_{13} = G^1 \circ (1 \mapsto 1, 3 \mapsto 2) =
\{ (1 \mapsto x, 3 \mapsto y): xy \in G^1) \}$
and $(G^H)^2_{23} = G^2$ (for either choice of $\tT_{23}$
if arcs of colour $2$ always come in opposite pairs).
All other $(G^H)^d_B$ are undefined.
If $(\Phi,G^H)$ is $(\oO,h)$-extendable
then for any sets $S_i \sub T_i$, $i \in [3]$
of size at most $h$ and an injection 
$\phi: S:=\bigcup_{i=1}^3 S_i \to [n]$
there are at least $\oO n^{|T \sm S|}$
injections $\phi^+: T:=\bigcup_{i=1}^3 T_i \to [n]$
extending $\phi$ such that for any $1 \le i<j \le 3$,
$x_i \in T_i$, $x_j \in T_j$ with $x_ix_j \not\sub S$ 
we have $(i \mapsto \phi^+(x_i), j \mapsto \phi^+(x_j))
\in (G^H)^d_{ij}$, i.e.\ $\phi^+(x_i)\phi^+(x_j) \in G^d$,
where $d=2$ if $ij=23$ or $d=1$ otherwise.

This is roughly equivalent to the following property:
for any disjoint $A,B \sub [n]$ of size at most $h$
there are at least $\oO n$ vertices $c$
such that $ca \in G^2$ for all $a \in A$
and $bc \in G^1$ for all $b \in B$,
and at least $\oO n$ vertices $c$
such that $ca \in G^1$ for all $a \in A \cup B$.

\medskip

\item
Recall the example of the two-coloured cyclic directed 
triangle: $q=3$, $r=2$, $t=1$, $D=2$, $\mc{H}=\{H\}$, 
$H^1 = \{ (1 \mapsto 1, 2 \mapsto 2),
(1 \mapsto 2, 2 \mapsto 3) \}$,
$H^2 = \{ (1 \mapsto 3, 2 \mapsto 1)\}$,
$\LL^1 = \LL^2 = \{id\}$. 
Let $\Phi$ be the complete $[3]$-complex on $[n]$
and $G \in (\mb{N}^2)^{\Phi_2}$.
The $2$-divisibility condition is trivial.
As $G(\es)^*=(|G^1|,|G^2|)$ and $H(\es)^*=(2,1)$
the $0$-divisibility condition is $|G^1|=2|G^2|$.
For $\psi \in I^1_n$, $x=Im(\psi) \in [n]$ we have
$G(\psi)^* = (|G^1(1 \mapsto x)|, |G^1(2 \mapsto x)|,
|G^2(1 \mapsto x)|, |G^2(2 \mapsto x)| )
= ( d_{G^1}^+(x), d_{G^1}^-(x), d_{G^2}^+(x), d_{G^2}^-(x) )$.
We have $H(1 \mapsto 1)^* = (1,0,0,1)$,
$H(1 \mapsto 2)^* = (1,1,0,0)$ and
$H(1 \mapsto 3)^* = (0,1,1,0)$, which generate
$H\sgen{1} = \{ \bm{v} \in \mb{Z}^4:
v_1 + v_3 = v_2 + v_4 \}$, 
so the $1$-divisibility condition is 
$d_{G^1}^+(x) + d_{G^2}^+(x) = d_{G^1}^-(x) + d_{G^2}^-(x)$,
i.e.\ the degree regularity condition 
$d_G^+(x)=d_G^-(x)$ needed for decomposition
into cyclic triangles ignoring the colours.

As for extendability, we have $(G^H)^1_{12} = G^1$,
$(G^H)^1_{23} = \{ (2 \mapsto x, 3 \mapsto y): xy \in G^1) \}$
and $(G^H)^2_{13} = \{ (3 \mapsto x, 1 \mapsto y): xy \in G^2) \}$.
If $(\Phi,G^H)$ is $(\oO,h)$-extendable then for any $S_i$, $T_i$ 
and $\phi$ as above there are at least $\oO n^{|T \sm S|}$ 
extensions $\phi^+$ such that for any $1 \le i<j \le 3$,
$x_i \in T_i$, $x_j \in T_j$ with $x_ix_j \not\sub S$ 
we have $\phi^+(x_i)\phi^+(x_j) \in G^1$ if $ij \ne 13$
or $\phi^+(x_3)\phi^+(x_1) \in G^2$ if $ij=13$.

This is roughly equivalent to:
for any disjoint $A,B \sub [n]$ of size at most $h$ there are 

(1) at least $\oO n$ vertices $c$
such that $ca \in G^1$ for all $a \in A$
and $bc \in G^1$ for all $b \in B$,

(2) at least $\oO n$ vertices $c$
such that $ca \in G^1$ for all $a \in A$
and $bc \in G^2$ for all $b \in B$, and

(3) at least $\oO n$ vertices $c$
such that $ca \in G^2$ for all $a \in A$
and $bc \in G^1$ for all $b \in B$.

\medskip

\item
Recall the example of the cyclic triangle
$1 \to 2 \to 3 \to 1$ with vertex partition 
$\mc{P}=(\{1,2\},\{3\})$: we have
$q=3$, $r=2$, $t=2$, $D=3$, $\mc{H}=\{H\}$, 
$H^1 = \{ (1 \mapsto 1, 2 \mapsto 2) \}$,
$H^2 = \{ (1 \mapsto 1, 2 \mapsto 3) \}$, 
$H^3 = \{ (1 \mapsto 2, 2 \mapsto 3) \}$,
$\iB^1=(2,0)$, $R((2,0))=([2],\es)$,
$\iB^2=\iB^3=(1,1)$, $R((1,1))=(\{1\},\{2\})$,
$\LL^1=\LL^2=\LL^3=\{id\}$. 
Let $\Phi$ be a complete $\mc{P}$-partite 
$[3]$-complex and $G \in (\mb{N}^3)^{\Phi_2}$.
Note that $P'_1 = \Phi^\circ_{\{1\}} = \Phi^\circ_{\{2\}}$
and $P'_2 = \Phi^\circ_{\{3\}}$.
The $2$-divisibility condition is that arcs of $G$
must respect the partition according to their colour,
i.e.\ if $G^1_\tT \ne 0$ then $Im(\tT) \sub P'_1$ and
if $G^2_\tT \ne 0$ or $G^3_\tT \ne 0$
then $\tT(1) \in P'_1$ and $\tT(2) \in P'_2$.
As $G(\es)^*=(|G^1|,|G^2|,|G^3|)$ and $H(\es)^*=(1,1,1)$
the $0$-divisibility condition is $|G^1|=|G^2|=|G^3|$,
i.e.\ in the uncoloured interpretation
we have equal numbers of arcs (1) within $P'_1$,
(2) from $P'_1$ to $P'_2$, and (3) from $P'_2$ to $P'_1$.

Now consider the $1$-divisibility conditions.
Let $G'$ denote the arcs between $P'_1$ and $P'_2$
according to the uncoloured interpretation,
where arcs from $P'_1$ to $P'_2$ correspond to $G^3$
and arcs from $P'_2$ to $P'_1$ correspond to $G^2$.
Let $\psi \in I^1_n$ and $x=Im(\psi)$.
Suppose first that $x \in P'_1$. Then
$G(\psi)^* = (|G^1(1 \mapsto x)|, |G^1(2 \mapsto x)|,
|G^2(1 \mapsto x)|, |G^2(2 \mapsto x)|,
|G^3(1 \mapsto x)|, |G^3(2 \mapsto x)| )
= ( d_{G^1}^+(x), d_{G^1}^-(x), 
d_{G'}^-(x), 0,  d_{G'}^+(x), 0 )$.
As $H(1 \mapsto 1)^* = (1,0,1,0,0,0)$
and $H(1 \mapsto 2)^* = (0,1,0,0,1,0)$,
we obtain the conditions 
$d_{G^1}^+(x)=d_{G'}^-(x)$ and
$d_{G^1}^-(x)=d_{G'}^+(x)$ for all $x \in P'_1$.
Now suppose $x \in P'_2$. We have
$G(\psi)^* = ( 0, 0, 0, d_{G'}^+(x), 0, d_{G'}^-(x))$
and $H(1 \mapsto 3)^* = (0,0,0,1,0,1)$,
so we need $d_{G'}^+(x)=d_{G'}^-(x)$ for all $x \in P'_2$.

As for extendability, we have $(G^H)^1_{12} = G^1$,
$(G^H)^2_{13} = \{ (1 \mapsto x, 3 \mapsto y): xy \in G^2) \}$
and $(G^H)^3_{23} = \{ (2 \mapsto x, 3 \mapsto y): xy \in G^3) \}$.
If $(\Phi,G^H)$ is $(\oO,h)$-extendable then for any $S_i$, $T_i$ 
and $\phi$ as above with $\phi(S_1), \phi(S_2) \sub P'_1$ 
and $\phi(S_3) \sub P'_2$ there are at least 
$\oO |P'_1|^{|T_1 \sm S_1| + |T_2 \sm S_2|} |P'_2|^{|T_3 \sm S_3|}$ 
extensions $\phi^+$ such that for any $1 \le i<j \le 3$,
$x_i \in T_i$, $x_j \in T_j$ with $x_ix_j \not\sub S$ 
we have $\phi^+(x_i)\phi^+(x_j) \in G^{d_{ij}}$,
where $d_{12}=1$, $d_{13}=2$, $d_{23}=3$.
This is roughly equivalent to:

(1) for any disjoint $A,B \sub P'_1$ of size at most $h$
there are at least $\oO |P'_2|$ vertices $c \in P'_2$
such that $ca \in G'$ for all $a \in A$
and $bc \in G'$ for all $b \in B$, and

(2) for any disjoint $A \sub P'_1$, $B \sub P'_2$ of size 
at most $h$ there are at least $\oO |P'_1|$ vertices 
$c \in P'_1$ such that $ca \in G'$ for all $a \in A$
and $bc \in G^1$ for all $b \in B$,
and at least $\oO |P'_1|$ vertices 
$c \in P'_1$ such that $ca \in G^1$ for all $a \in A$
and $bc \in G'$ for all $b \in B$.
\end{enumerate}

\medskip

Similarly to Definition \ref{def:diext=},
we have the following general rough equivalence: 
if $(\Phi,G^H)$ is $(\oO,h)$-extendable
then $(\Phi,G)$ is $(\oO,h,H)$-vertex-extendable
(as in the next definition),
and conversely, if $G$ is $(\oO,h,H)$-vertex-extendable
then $(\Phi,G^H)$ is $(\oO^{qh},h)$-extendable.

\begin{defn} \label{def:masterext=}
With notation as in Definition \ref{def:master},
we say $(\Phi,G)$ is $(\oO,h,H)$-vertex-extendable
if for any $x \in [q]$ and disjoint sets $A_i$, 
$i \in [q] \sm \{x\}$ of size at most $h$ such that
$(i \mapsto v_i: i \in [q]\sm \{x\}) \in \Phi$
whenever each $v_i \in A_i$, there are at least $\oO n$ 
vertices $v \in \Phi^\circ_x$ such that
\begin{enumerate}
\item $(i \mapsto v_i: i \in [q]) \in \Phi$
whenever $v_x=v$ and $v_i \in A_i$ for each $i \ne x$,
\item for each $d \in [D]$ and 
arc $\tT$ of $H^d$ with $x \in Im(\tT)$,
we have all arcs $(i \mapsto v_i: i \in [r])$ in $G^d$
where $v_j=v$ for $j=\tT^{-1}(x)$
and $v_i \in A_{\tT(i)}$ for all $i \ne j$.
\end{enumerate}
\end{defn}

The main theorem of the section provides 
the above general setting with our usual conclusion 
(divisibility, regularity and extendability
suffice for the existence of decompositions).

\begin{theo} \label{master}
With notation as in Definition \ref{def:master},
suppose all $n_1/h \le |P'_i| \le n_1$ with $n_1 > n_0(q,D)$,
that $G$ is $\mc{H}$-divisible and
$(\mc{H},c,\oO)$-regular in $\Phi$,
and all $(\Phi,G^H)$ are $(\oO,h)$-extendable,
where $n_1^{-\dD}<\oO<\oO_0(q,D)$ and $c=\oO^{h^{20}}$.
Then $G$ has an $\mc{H}$-decomposition in $\Phi$.
\end{theo}

\nib{Proof.} 
For $\psi \in \Phi_{[r]}$ we let $\psi^*$ be the set of all 
$\psi \circ \pi^{-1}$ where $\pi:[r] \to [q]$ is order-preserving
and $i_{\mc{P}}(\pi)=i_{\mc{P}'}(\psi)$.
Similarly, for $\tT \in H_r$ we let $\tT^*$ be the set of all 
$\tT \circ \pi^{-1}$ where $\pi:[r] \to [q]$ is order-preserving
and $i_{\mc{P}}(\pi)=i_{\mc{P}}(\tT)$.
Let $G^* = \sum_{\psi \in \Phi_{[r]}} G_\psi \psi^*$ 
and $\mc{H}=\{H^*: H \in \mc{H}\}$ with 
each $(H^*)^d=(H^d)^*=\{\tT^*: \tT \in H^d\}$.
Let $\mc{A}=\{A^H: H \in \mc{H}\}$ with each $A^H=\Ss^\le$ and
$\gG \in \mb{Z}^{\mc{A}_r}$ where each $\gG_\tT$ is $e_d$
if $\tT \in H^{d*}$ for some $H \in \mc{H}$, $d \in [D]$,
otherwise zero. Then a $\gG(\Phi)$-decomposition of $G^*$
is equivalent to an $\mc{H}^*$-decomposition of $G^*$, 
and so, we claim, to an $\mc{H}$-decomposition of $G$.

For the latter equivalence, similarly to Lemma
\ref{equivdi*}, we need to show for any 
$H \in \mc{H}$, $d \in [D]$, $\phi \in \Phi_q$
that $\psi \in \phi H^d_r$ iff $\psi^* \sub \phi H^{d*}_r$.
To see this, write $\psi = \phi\tT$, where $\tT \in H^d_r$
and let $\iB=i_{\mc{P}'}(\psi)=i_{\mc{P}}(\tT)$.
For any $\psi' \in \psi^*$ we can write $\psi'=\psi\pi^{-1}$
where $\pi:[r] \to [q]$ is order-preserving
with $i_{\mc{P}}(\pi) = \iB$,
so $\psi'=\phi\tT'$ with $\tT'=\tT\pi^{-1} \in \tT^*$.
Thus $\psi \in \phi H^d_r$ implies $\psi^* \sub \phi H^{d*}_r$.
The converse is similar, so the claimed equivalence holds
(and also for integral decompositions).

Next we claim that $\gG$ is elementary. To see this, 
we describe the type vectors $\gG^\tT$ 
for $\tT \in A^H_B$, $B \in Q$. If $\gG^\tT \ne 0$ 
then we can write $\tT = \tT_0 \tau_0 \pi_0^{-1}$ 
with $\tT_0 \in H$, $\tau_0 \in S_r$ 
and $\pi_0 \in Bij([r],B)$ order-preserving.
Say $\tT_0 \in H^d$. As $\mc{H}$ is $(\mc{P},\LL)$-canonical, 
for $\tT' \in Bij([r],Im(\tT_0))$ we have 
$\tT' \in H^d$ iff $\tT^{-1} \tT' \in \LL^d$.
Fix a set $X^d$ of representatives for the right cosets
of $\LL^d$ in $S_r$. Then we have a unique expression
$\tT = \tT_0 \tau_0 \pi_0^{-1}$ with $\tT_0 \in H^d$ 
and $\tau_0 \in X^d$. For any $\sS \in \Ss^B$ we have 
$\gG^\tT_\sS = \gG_{\tT\sS} \in \{0,e_d\}$ equal to $e_d$ 
iff $\sS = \pi_0 (\lL\tau_0)^{-1} \pi^{-1}$ where $\lL \in \LL^d$
and $\pi:[r] \to [q]$ is order-preserving
with $i_{\mc{P}}(\pi)=\iB:=i_{\mc{P}}(B)$.
Thus, besides the $0$ type, for each $B \in Q$
and $d \in [D]$ with $\iB^d=\iB$ we have
$|X^d|=r!/|\LL^d|$ types $(t^{\tau_0}: \tau_0 \in X^d)$
describing all generalised arcs with any given image.
Given $d$, the supports of the $t^{\tau_0}$ for $\tau_0 \in X^d$ 
are mutually disjoint, so $\gG$ is elementary, as claimed.

As $G$ is $(\mc{H},c,\oO)$-regular in $\Phi$ we have
$y^H_\phi \in [\oO n^{r-q},\oO^{-1} n^{r-q}]$ for each 
$H \in \mc{H}$, $\phi \in \Phi_q$ with $\phi H \le G$
(equivalently, $\phi H^* \le G^*$) so that 
$\sum \{ y^H_\phi \phi H \} = (1 \pm c)G$.
(equivalently, $\sum \{ y^H_\phi \phi H^* \} = (1 \pm c)G^*$).
We identify any such $\phi H^* \le G^*$ with
$\gG(\phi) \le_\gG G^*$ (regarding $\phi \in A^H(\Phi)$),
so $\phi \in \mc{A}(\Phi,G^*)$. 
Let  $y_\phi = y^H_\phi$ for $\phi \in A^H(\Phi)$.
For any $B \in Q$, $\psi \in \Phi_B$, 
$d \in [D]$ with $\iB^d=\iB:=i_{\mc{P}'}(\psi)$
and $0 \ne t \in T_B$, say with $t$ supported on the set 
of all $(\lL\tau)^{-1} \pi^{-1}$ where $\lL \in \LL^d$ 
and $\pi:[r] \to [q]$ is order-preserving 
with $i_{\mc{P}}(\pi)=\iB$, we have
$\pl^t y_\psi = \sum \{ y_\phi : t_\phi(\psi)=t \} 
= \sum \{ y^H_\phi: \psi\tau \in \phi H^{d*}, H \in \mc{H} \}
= (1 \pm c) G^{d*}_{\psi\tau} = (1 \pm c)(G^*)^t_\psi$,
so $G^*$ is $(\gG,c,\oO)$-regular.

Next we consider extendability. Fix $H \in \mc{H}$. We have
$\gG[G^*]^H = \{ \psi \in A^H(\Phi)_r : \gG(\psi) \le_\gG G^* \}$,
so $\psi \in \Phi_B$ is in $\gG[G^*]^H$ iff 
(a) no arc in $H$ has image $B$, or
(b) $\psi \tT_B \in G^d$ (i.e.\ $G^d_{\psi \tT_B}>0$) for some 
(equivalently, all) $\tT_B \in H^d_r$ with $Im(\tT_B)=B$. 
Let $E=(J,F,\phi)$ be any $\Phi$-extension of rank $h$
and $J' \sub J_r \sm J[F]$. Let $J'' = (J^d: d \in [D])$ 
with each $J^d = \bigcup \{ J'_B: \tT_B \in H^d_r \}$.
As $(\Phi,G^H)$ is $(\oO,h)$-extendable
we have $X_{E,J''}(\Phi,G^H) > \oO n^{v_E}$.
Consider any $\phi^+ \in X_{E,J''}(\Phi,G^H)$.
For any $\psi \in J^d_B$, $d \in [D]$ we have 
$\phi^+\psi  \in (G^H)^d_B$, 
so $\phi^+\psi \tT_B  \in G^d$, 
so $\phi^+\psi \in \gG[G^*]^H$.
Thus $\phi^+ \in X_{E,J'}(\Phi,\gG[G^*]^H)$,
so $(\Phi,\gG[G^*])$ is $(\oO,h)$-extendable.

To deduce the theorem from Theorem \ref{main},
it remains to show for any orbit $O \in \Phi/\Ss$
that $((G^*)^\sharp)^O \in \sgen{\gG^\sharp[O]}$.
Fix $\psi \in O \in \Phi_i/\Ss$. 
Let $\iB'=i_{\mc{P}'}(\psi)$ and 
$I' = \{ \tT \in I^i_q: i_{\mc{P}}(\tT)=\iB'\}$.
Write $\psi = \psi_0 \pi_0^{-1}$ with $\psi_0 \in I^i_n$
and $\pi_0:[i] \to Dom(\psi)$ order-preserving.
As $G$ is $\mc{H}$-divisible, there is 
$n \in \mb{Z}^{\mc{H} \times I'}$ with
$G(\psi_0)^* = \sum_{H,\tT} n_{H\tT} H(\tT)^*$,
i.e.\ $|G^d(\psi_0 \pi^{-1})|
= \sum_{H,\tT} n_{H\tT} |H^d(\tT \pi^{-1})|$
for all $d \in [D]$ and $\pi \in I^i_r$.
Writing $\phi = \psi_0 \tT^{-1} \in A^H(\Phi)$, we claim 
$((G^*)^\sharp)^O = \sum_{H,\tT} n_{H,\tT} \gG^\sharp(\phi)$.
To see this, fix $\psi\sS \in O$, 
$B \in Q$ and let $\iB=i_{\mc{P}}(B)$.
We need to show for any $d \in [D]$ with $\iB^d=\iB$ 
that $|G^{d*}_B(\psi\sS)| = \sum_{H,\tT} 
n_{H,\tT} |H^{d*}_B(\tT\pi_0^{-1}\sS)|$, 
i.e.\ $|G^d(\psi_0\pi^{-1})|
= \sum_{H,\tT} n_{H\tT} |H^d(\tT\pi^{-1})|$,
where $\pi^{-1} = \pi_0^{-1} \sS \pi_B$
with $\pi_B:[r] \to B$ order-preserving;
this is a case of the previous identity.
\qed

\section{Perspectives} \label{sec:end}

The existence of designs established in \cite{Kexist}
has seen several subsequent applications,
some of which are particularly instructive 
as they require not only the existence but also that
designs can be `almost entirely random', in that the semi-random (nibble) 
construction of approximate designs by R\"odl \cite{R} can be completed 
to an actual design by an absorption process 
(Randomised Algebraic Construction in \cite{Kexist} 
or Iterative Absorption in \cite{GKLO}).
In this vein, we mention the proof by Kwan \cite{Kw} that
almost all Steiner triple systems have perfect matchings,
results on discrepancy of high-dimensional permutations
by Linial and Luria \cite{LL2}, and the existence of bounded degree 
coboundary expanders of every dimension
by Lubotzky, Luria and Rosenthal \cite{LLR}.
These results suggest that the new results in \cite{K2} may create
more fruitful connections with the theory of high-dimensional expanders
and other topics in high-dimensional combinatorics.

In Design Theory, the most fundamental problems that remain open
are those concerning designs with large block sizes. Here we recall
from the introduction the Prime Power Conjecture on projective planes,
where we know that the divisibility conditions do not always suffice;
the conjecture seems to reflect a philosophy that a combinatorial
description of a sufficient rich structure 
somehow implies an algebraic characterisation.
On the other hand, a conjecture that reflects the opposite philosophy
is that Hadamard matrices (see \cite{had}) of order $n$ should
exist whenever the trivially necessary conditions are satisfied
(i.e.\ $n$ is $1$, $2$ or divisible by $4$). It is not clear how the
methods of \cite{GKLO,GKLO2,Kexist,K2} could apply to such problems,
where a more fruitful direction may be the development of the
approach of \cite{KLP}, which can allow for large block sizes.
There are also many well-known open problems in Design Theory that
do not involve large block sizes, and so may be more approachable
by absorption techniques. Here we mention Ryser's Conjecture \cite{ryser}
that every Latin square of odd order should have a transversal;
equivalently, any triangle decomposition of $K_3(n)$ for $n$ odd
should contain a triangle factor (perfect matching of triangles).

In Combinatorics, there are several natural directions in which
one may seek to generalise the existence of various types of design,
from extremal and/or probabilistic perspectives. A basic class of 
extremal questions is to determine the minimum degree threshold
(which has various possible definitions) for decompositions
(see e.g.\ \cite{GKLO2,M}). Natural probabilistic directions
are thresholds for the existence of certain designs 
in random hypergraphs (e.g.\ Steiner Triple Systems in $G^3(n,p)$)
or a theory of Random Designs analogous
to the rich theory of Random Graphs.

\medskip

\nib{Acknowledgement.} 
I would like to thank an anonymous referee for very detailed and helpful comments on the presentation of this paper.

\end{document}